\documentclass[final,a4paper,leqno,11pt]{article}
\usepackage{graphicx}
\usepackage{mathptmx}
\usepackage{amssymb,amsmath,amsthm, amsfonts}
\usepackage{a4wide}
\usepackage{color}

\DeclareMathOperator{\rank}{rank}
\newcommand{\R}{\mathbb{R}}
\newcommand{\N}{\mathbb{N}}
\newcommand{\norm}[1]{\|#1\|}

\newcommand{\dist}[1]{{\rm d}(#1)}
\newcommand{\B}{{\cal B}}

\newcommand{\mv}{\,\vert\, }

\newcommand{\Lag}{{\cal L}}
\newcommand{\setto}[1]{\mathop{\to}\limits^#1}
\newcommand{\skalp}[1]{\langle #1\rangle}
\newcommand{\xb}{\bar x}
\newcommand{\AT}[2]{{\textstyle{#1\atop#2}}}
\newcommand{\xba}{{\bar x^\ast}}
\newcommand{\oo}{o}
\newcommand{\I}{{\cal I}}
\newcommand{\J}{{\cal C}}
\newcommand{\A}{{\cal A}}

\newcommand{\E}{{\cal E}}
\newcommand{\Tlin}{T_\Gamma^{\rm lin}}
\newcommand{\co}{{\rm conv\,}}
\newcommand{\Gr}{{\rm gph\,}}
\newcommand{\dom}{{\rm dom\,}}
\newcommand{\lip}{{\rm lip\,}}
\newcommand{\argmin}{{\rm argmin\,}}
\newcommand{\argmax}{\mathop{\rm arg\,max}\limits}

\def\ph{\hat{p}}

\def\xb{\bar{x}}

\def\A{\mathcal{A}}

\def\N{\mathcal{N}}

\def\R{\mathcal{R}}

\def\argmin{\mathop{\rm argmin}}
\def\argmax{\mathop{\rm argmax}}

\def\dist{\mathop{\rm dist}}
\def\dom{\mathop{{\rm dom}}}

\def\Gr{\mbox{\rm gph}\,}

\def\1B{{\bf 1}}
\def\R{\mathbb{R}}

\newcommand\be{\begin{equation}}
\newcommand\ee{\end{equation}}
\newcommand\ba{\begin{array}}
\newcommand\ea{\end{array}}
\newcommand{\bean}{\begin{eqnarray*}}
\newcommand{\eean}{\end{eqnarray*}}

\def\tilt{\mathop{\rm tilt}}

\def\disp{\displaystyle}

\def\la{\langle}
\def\ra{\rangle}

\def\st{\tilde{s}}

\def\y0{\bar Y_{0}}

\def\sce{\setcounter{equation}{0}}
\def\argmin{\mathop{{\rm argmin}}}

\def\sce{\setcounter{equation}{0}}

\def\sce{\setcounter{equation}{0}}
\def\disp{\displaystyle}

\def\tto{\;{\lower 1pt \hbox{$\rightarrow$}}\kern -10pt
\hbox{\raise 2pt \hbox{$\rightarrow$}}\;}
\def\Hat{\widehat}
\def\hat{\widehat}
\def\Tilde{\widetilde}
\def\tilde{\widetilde}
\def\Bar{\overline}
\def\ra{\rangle}
\def\la{\langle}

\def\B{\mathbb B}

\def\N{\mathcal{N}}

\def\h{\hfill\Box}
\def\R{I\!\!R}

\def\ox{\bar{x}}

\def\oz{\bar{z}}
\def\ov{\bar{v}}

\def\ow{\bar{w}}

\def\co{\mbox{\rm co}\,}

\def\gph{\mbox{\rm gph}\,}
\def\dist{\mbox{\rm dist}\,}

\def\dom{\mbox{\rm dom}\,}

\def\lip{\mbox{\rm lip}\,}

\def\rge{\mbox{\rm rge}\,}

\def\h{\hfill\triangle}
\def\dn{\downarrow}
\def\O{\Omega}

\def\ph{\varphi}
\def\emp{\emptyset}
\def\st{\stackrel}
\def\oR{\Bar{\R}}

\def\lm{\lambda}
\def\gg{\gamma}
\def\dd{\delta}

\def\bb{\beta}
\def\kk{\kappa}
\def\Th{\Theta}

\def\vt{\vartheta}
\def\R{\mathbb R}
\def\N{\mathbb N}
\def\sce{\setcounter{equation}{0}}

\date{}
\title{Complete Characterizations of Tilt Stability in Nonlinear Programming under Weakest Qualification Conditions}
\author{HELMUT GFRERER\footnote{Institute of Computational Mathematics, Johannes Kepler University Linz, A-4040 Linz, Austria (helmut.gfrerer@jku.at)} $\;$ and $\;$ BORIS S. MORDUKHOVICH\footnote{Department of Mathematics, Wayne State University, Detroit, MI 48202, USA (boris@math.wayne.edu).}}
\begin{document}
\newtheorem{Theorem}{Theorem}[section]
\newtheorem{Proposition}[Theorem]{Proposition}
\newtheorem{Remark}[Theorem]{Remark}
\newtheorem{Lemma}[Theorem]{Lemma}
\newtheorem{Corollary}[Theorem]{Corollary}
\newtheorem{Definition}[Theorem]{Definition}
\newtheorem{Example}[Theorem]{Example}
\renewcommand{\theequation}{{\thesection}.\arabic{equation}}
\maketitle

\small
{\bf Abstract.} This paper is devoted to the study of tilt stability of local minimizers for classical nonlinear programs with equality and inequality constraints in finite dimensions described by twice continuously differentiable functions. The importance of tilt stability has been well recognized from both theoretical and numerical perspectives of optimization, and this area of research has drawn much attention in the literature, especially in recent years. Based on advanced techniques of variational analysis and generalized differentiation, we derive here complete pointbased second-order characterizations of tilt-stable minimizers entirely in terms of the initial program data under the new qualification conditions, which are the weakest ones for the study of tilt stability.
\vspace*{0.05in}

\vspace*{0.05in}

\vspace*{0.05in}

{\bf Abbreviated title.} Tilt stability in nonlinear programming

\normalsize

\section{Introduction}\sce

The notion of {\em tilt-stable local minimizers} was introduced by Poliquin and Rockafellar \cite{PolRo98} for problems of unconstrained optimization with a general extended-real-valued objective function, which implicitly incorporates constraints via the indicator function of the feasible region. Motivated by the justification of convergence properties, stopping criteria, and robustness of numerical algorithms, the authors of \cite{PolRo98} suggested to study and characterize not just arbitrary local minimizers but those which behave nicely with respect to linear perturbations tilted the objective function in one direction or another; namely, minimizers that remain locally unique and Lipschitz continuous under small perturbations of the aforementioned type. Tilt stability has attracted strong attention in the literature, particularly in recent years; see, e.g., \cite{BS00,DrLe13,DMN14,EbWe12,LPR00,LeZh13,MoNg14,MoOut13,MoRo12,ZN14} and the references therein.

In \cite{PolRo98}, Poliquin and Rockafellar obtained a characterization of tilt-stable local minimizers for a large class of prox-regular extended-real-valued functions via the positive-definiteness of their {\em second-order subdifferential/generalized Hessian} in the sense of Mordukhovich \cite{Mo92}; see Section~2. Based on this result and the newly developed second-order calculus rules, Mordukhovich and Rockafellar \cite{MoRo12} derived a characterization of tilt-stable local minimizers for nonlinear programs (NLPs) with $C^2$-smooth data assuming the {\em linear independence constraint qualification} (LICQ). Under this nondegeneracy assumption, the characterization of tilt stability was expressed in \cite{MoRo12} via Robinson's {\em strong second-order sufficient condition} (SSOSC) \cite{Rob80} formulated entirely {\em at} the local minimizer in question; such conditions are called {\em pointbased} (known also as {\em pointwise}) in what follows. They are surely much more preferable for applications than the {\em neighborhood} conditions discussed below.

In the further lines of research, Mordukhovich and Nghia \cite{MoNg14} introduced the notion of tilt-stable local minimizers with {\em modulus $\kk>0$} for an extended-real-valued objective function and derived, by developing a new {\em dual-space} approach to tilt stability, a characterization of such minimizers in terms of the so-called {\em combined second-order subdifferential} (see Section~2) via a strong positive-definiteness condition involving $\kk$. It is shown in \cite{MoNg14} that the obtained characterization reduces to the one in \cite{PolRo98} when the modulus $\kk$ is not an issue, and also that the aforementioned result of \cite{MoNg14} leads to a new characterization of tilt-stable local minimizers for NLPs with $C^2$-smooth inequality constraints without imposing LICQ. Namely, the LICQ assumption was weakened in \cite{MoNg14} by the simultaneous fulfillment of the {\em Mangasarian-Fromovitz constraint qualification} (MFCQ) and the {\em constant rank constraint qualification} (CRCQ) conditions, while the characterization of tilt-stable local minimizers was given in this setting via the new {\em uniform second-order sufficient condition} (USOSC). The new USOSC is shown in \cite{MoNg14} to be strictly weaker than SSOSC, while being reduced to the latter under the validity of LICQ. However, in contrast to LICQ, MFCQ, and SSOSC, the formulations of CRCQ and USOSC are not pointbased depending on points in a {\em neighborhood} of the reference local minimizer.

As demonstrated by simple examples (see Section~8), the combination of MFCQ and CRCQ constitutes a setting, which is not fully satisfactory for the study of tilt-stable minimizers and may exclude from consideration important situations when tilt-stable minimizers exist and can be recognized. Furthermore, the obtained USOSC characterization in \cite{MoNg14} is a neighborhood condition but not a pointbased one. On the other hand, the results below show that under the MFCQ assumption {\em alone} a pointbased second-order characterization of tilt stability is {\em not possible}, which means that there are two NLPs with the same derivatives up to the second order at the reference point satisfying MFCQ but such that one problem admits a tilt-stable minimizer at this point while the other one doesn't.

To go forward in this paper, we dispense with MFCQ and also with CRCQ by replacing them, in the general case of both $C^2$-smooth {\em inequality} and {\em equality} constraints in NLPs, with another pair of constraint qualifications such that the simultaneous fulfillment of these conditions is {\em strictly weaker} than the validity of {\em each} of the conditions MFCQ and CRCQ and thus of their combination.  The first of these new assumptions/qualification conditions, called the {\em metric subregularity constraint qualification} (MSCQ), weakens the property of metric regularity for the constraint NLP mapping {\em around} the reference minimizer (the latter property is known to be equivalent to MFCQ of this mapping at the point in question) by its metric {\em sub}regularity {\em at} this point. This assumption has been recently employed in the papers by Gfrerer and Outrata \cite{GfrOut14a,GfrOut14b} for evaluating generalized derivatives of the normal cone mapping to inequality systems. An effective pointbased condition for its validity in terms of the first-order and second-order derivatives of the inequality constraint functions was introduced earlier by Gfrerer \cite{Gfr11} as the {\em second-order sufficient condition for metric subregularity} (SOSCMS).

The other qualification condition coupled here with (and independent of) MSCQ appears for the first time in this paper under the name of the {\em bounded extreme point property} (BEPP) of the NLP constraint system. Although it is formulated at the reference solution, neighborhood points are used in the definition as well. Being much weaker than MFCQ and CRCQ, this new qualification condition is also implied by the aforementioned SOSCMS property, which is completely pointbased.

Involving only the weakest qualification conditions MSCQ and BEPP, the main results of the paper provide {\em pointbased second-order characterizations} of tilt-stable minimizers for general NLPs with $C^2$-smooth equality and inequality constraints entirely in terms of their initial data. These characterizations are given as follows: first we derive pointbased {\em sufficient} conditions for tilt-stable minimizers in full generality and then show that they are also {\em necessary} for tilt stability under some additional assumptions. The major sufficient condition for tilt-stable minimizers is formulated via the {\em strong positive-definiteness} (depending on modulus $\kappa>0$ of tilt stability and being {\em much weaker} than SSOSC) of the Hessian of the Lagrange function by using a rather narrow subset of {\em extreme points} of the collection of those Lagrange multipliers, which are solutions to a certain linear program associated with {\em critical directions}. The {\em necessity} of the second-order conditions is justified under different additional assumptions: either {\em nondegeneracy} in critical directions, or the so-called {\em 2-regularity}, or {\em CRCQ}. Furthermore, we show that even without these extra assumptions the aforementioned {\em sufficient} conditions become {\em necessary} in a slightly modified problem with the same cost function and constraints reducing to the original ones together with their first and second derivatives at the reference minimizers. We also provide the {\em quantitative evaluation} (estimates and exact formulas) of the corresponding moduli. All of this allows us to conclude that the obtained second-order pointbased characterizations of tilt-stable minimizers in NLPs are {\em complete}.\vspace*{0.05in}

The rest of the paper is organized as follows. Section~2 presents some basic constructions and properties of variational analysis and generalized differentiation widely used in the main body of the paper. We also formulate here the notion of {\em tilt-stable minimizers} and its {\em neighborhood characterization} in the unconstrained extended-real-valued format of optimization.

Section~3 is devoted to {\em qualification conditions} for NLPs with smooth inequality and equality constraints. We define here the main MSCQ and BEPP conditions and establish their connections with constraint qualifications well recognized in nonlinear programming and used in the paper. The next Section~4 involves from one side some {\em second-order analysis} to better understand both MSCQ and BEPP for NLPs with $C^2$-smooth data, while from the other side it demonstrates a significant role of these qualification conditions to get the desired properties of the indicator function of the constraint set and also to describe the {\em critical cone} to this set and the collection of Lagrange multipliers in {\em critical directions} needed for the subsequent second-order analysis of tilt-stable minimizers.

Section~5 contains important results on the precise calculation under the imposed weakest qualification conditions MSCQ and BEPP of some second-order generalized differential constructions for sets of feasible solutions to NLPs with $C^2$-smooth inequality and equality constraints via their initial data. These results, being certainly of their own interest, are used in this section for evaluating some second-order terms crucial for the subsequent characterizations of tilt-stable minimizers. Note that the obtained calculating formulas not only extend the corresponding results of \cite{GfrOut14a} to programs with equality constraints, but also replace a certain relaxed metric regularity condition used in \cite{GfrOut14a} for the case of inequalities by the weaker BEPP qualification condition developed in this paper.

Sections~6 and 7 are central in the paper containing the pointbased second-order conditions for tilt-stable minimizers described above; namely, sufficient conditions in Section~6 and rather general while different settings for their necessity and complete characterizations presented in Section~7.

The obtained results on tilt stability are discussed and illustrated by various examples in Section~8, where important features of the developed necessary and sufficient conditions and the imposed qualifications are revealed in comparison with known results in this direction while discussing also related numerical issues. The final Section~9 contains some open questions and topics for the future research.\vspace*{0.05in}

Our notation is basically standard in variational analysis and optimization; see, e.g., \cite{Mo06a,RoWe98}. Recall that $\B$ stands for the closed unit ball in the finite-dimensional space in question with the Euclidean norm $\|\cdot\|$ and the scalar product $\la\cdot,\cdot\ra$ between two vectors, $\B_r(x):=x+r\B$, the symbol $^*$ indicates in general a dual operation including the matrix transposition, the polar cone, etc., $\dist(x;\O)$ denotes the distance from the point $x$ to the set $\O$, the symbol $\nabla q(\ox)$ stands for the gradient of a scalar function and for the Jacobian matrix for a vector one at $\ox$, and finally we have $\N:=\{1,2,\ldots\}$.

\section{Preliminaries from Variational Analysis}\sce

Let $f:\R^n\to\oR:=\R\cup\{\infty\}$ be an extended-real-valued function, which is assumed to be {\em proper}, i.e., $\dom f:=\{x\in\R^n\mv f(x)<\infty\}\not=\emp$, and let $\xb\in\dom f$. The {\em regular subdifferential} (known also as the presubdifferential and as the Fr\'echet/viscosity subdifferential) is defined by
\begin{eqnarray}\label{rs}
\widehat\partial f(\xb):=\left\{v^\ast\in\R^n\Big|\;\liminf_{x\to\xb}\frac{f(x)-f(\xb)-\skalp{v^\ast,x-\xb}}{\norm{x-\xb}}\ge 0\right\}.
\end{eqnarray}
The {\em limiting subdifferential} (known also as the Mordukhovich/basic subdifferential) of $f$ at $\xb$ is
\begin{eqnarray}\label{ls}
\partial f(\xb):=\big\{v^\ast\in\R^n\big|\;\exists\;x_k\to\xb,\;v_k^\ast\to v^\ast\;\mbox{ with }\;f(x_k)\to f(\xb),\;v_k^\ast\in\Hat\partial f(x_k),\;k\in\N\big\}.
\end{eqnarray}

Both constructions \eqref{rs} and \eqref{ls} reduce to the subdifferential of convex analysis if $f$ is convex. For $C^1$-smooth functions the subdifferentials $\Hat\partial f(x)$ and $\partial f(x)$ consist only of the gradient $\nabla f(x)$.

A lower semicontinuous (l.s.c.) function $f:\R^n\to\oR$ is called {\em prox-regular} at $\xb\in\dom f$ for $\bar v^\ast\in\partial f(\xb)$ if there are reals $r,\epsilon>0$ such that for all $x,u\in\B_\epsilon(\xb)$ with $\vert f (u)-f(\xb)\vert\le\epsilon$ we have
\begin{eqnarray}\label{prox-reg}
f(x)\ge f(u)+\skalp{v^\ast,x-u}-\frac r2\norm{x-u}^2\;\mbox{ whenever }\;v^\ast\in\partial f(u)\cap\B_{\epsilon}(\bar v^\ast).
\end{eqnarray}
Such a function is said to be {\em subdifferentially continuous} at $\xb\in\dom f$ for $\bar v^\ast\in\partial f(\xb)$ if $f(x_k)\to f(\xb)$ for all sequences $x_k\to\xb$ and $v_k^\ast\to\bar v^\ast$ as $k\to\infty$ with $v_k^\ast\in\partial f(x_k)$, $k\in\N$.

In what follows we also need some concepts from variational geometry. Given a set $\Omega\subset\mathbb R^d$ and a point $\bar z\in\Omega$, define the (Bouligand-Severi) {\em tangent/contingent cone} to $\Omega$ at $\bar z$  by
\begin{equation}\label{EqTan}
T_\Omega(\bar z):=\big\{u\in\mathbb R^d\big|\;\exists\,t_k\downarrow 0,\;u_k\to u\;\mbox{ with }\;\bar z+t_k u_k\in\O\;\mbox{ for all }\;k\big\}.
\end{equation}
The (Fr\'{e}chet) {\em regular normal cone} to $\Omega$ at $\bar z\in\Omega$ can be equivalently defined either by
\begin{equation}\label{fn1}
\widehat N_\Omega(\bar
z):=\Big\{v^\ast\in\R^d\Big|\;\limsup_{z\stackrel{\Omega}{\to}\bar z}\frac{\skalp{v^\ast,z-\bar z}}{\|z-\bar z\|}\le 0\Big\},
\end{equation}
where $z\st{\O}{\to}\oz$ means that $z\to\oz$ with $z\in\O$, or as the dual/polar to the contingent cone \eqref{EqTan}, i.e., by
\begin{equation}\label{fn2}
\widehat N_\Omega(\bar z):=T_\Omega(\bar z)^*.
\end{equation}
For convenience, we put $\Hat N_\O(\oz):=\emp$ for $\oz\notin\O$. Further, the (Mordukhovich) {\em limiting/basic normal cone} to $\Omega$ at $\bar z\in\Omega$ is given by
\begin{eqnarray}\label{nc}
N_\O(\oz):=\big\{\ov^*\in\R^d\big|\,\exists\;
z_k\to\oz,\;v^*_k\to\ov^*\;\mbox{ with }\;v_k^\ast\in\widehat N_\Omega(z_k)\;\mbox{ for all }\;k\big\}.
\end{eqnarray}
Note that, in spite of (in fact due to) being nonconvex, the normal cone \eqref{nc} and the corresponding to it limiting subdifferential and coderivative constructions enjoy {\em full calculi}, which are based on variational/extremal principles of variational analysis; see, e.g., \cite{Mo06a,RoWe98} and the references therein.

It is easy to observe the following relationships:
\[
\widehat\partial\delta_\Omega(z)=\widehat N_\Omega(z)\;\mbox{ and }\;\partial\delta_\Omega(z)= N_\Omega(z)\;\mbox{ for all }\;z\in\O
\]
between the corresponding regular and limiting subdiffential/normal cone constructions defined above, where $\dd_\O(z)$ stands for the indicator function of the set $\O$ equal to $0$ if $z\in\Omega$ and to $\infty$ otherwise.

Considering next a set-valued (in particular, single-valued) mapping $\Psi:\R^d\rightrightarrows\R^s$, we associate with it the {\em domain} $\dom\Psi$ and the {\em graph} $\Gr\Psi$ by
\[
\dom\Psi:=\big\{z\in\R^d\big|\;\Psi(z)\ne\emp\big\}\;\mbox{ and }\;\Gr\Psi:=\big\{(z,w)\big|\;w\in\Psi(z)\big\}
\]
and define the following coderivative constructions. The {\em regular coderivative} $\widehat D^\ast\Psi(\bar z,\bar w):\R^s\rightrightarrows\R^d$ of $\Psi$ at $(\bar z,\bar w)\in\Gr\Psi$ is generated by the regular normal cone \eqref{fn1} via
\begin{equation}\label{EqRegCoderiv}
\widehat D^\ast\Psi(\bar z,\bar w)(v^\ast):=\big\{u^\ast\in\R^d\big|\;(u^\ast,-v^\ast)\in\widehat N_{{\rm{\small gph}}\,\Psi}(\bar z,\bar w)\big\},\quad v^\ast\in\R^s,
\end{equation}
and the {\em limiting coderivative} $\widehat D^\ast\Psi(\bar z,\bar w):\R^s\rightrightarrows\R^d$ of $\Psi$ at $(\bar z,\bar
w)\in\Gr\Psi$ is given by
\begin{equation}\label{EqLimCoderiv}
D^\ast\Psi(\bar z,\bar w)(v^\ast):=\big\{u^\ast\in\R^d\big|\;(u^\ast,-v^\ast)\in N_{{\rm{\small gph}}\,\Psi}(\bar z,\bar w)\big\},\quad v^\ast\in\R^s.
\end{equation}
If $\Psi$ is single-valued at $\bar z$, we drop $\bar w$ in the notation of \eqref{EqRegCoderiv} and \eqref{EqLimCoderiv}. If $\Psi$ is a single-valued mapping being smooth around $\oz$, then we have the equalities
$$
\widehat D^\ast\Psi(\bar z)(v^\ast)=D^\ast\Psi(\bar z)(v^\ast)=\big\{\nabla\Psi(\bar z)^*v^\ast\big\}\;\mbox{ for all }\;v^\ast\in\R^s.
$$

One of the striking advantages of the limiting coderivative \eqref{EqLimCoderiv} (besides full calculus) is the possibility to derive in its terms complete pointbased characterizations of some basic properties of {\em well-posedness} in nonlinear and variational analysis related to robust Lipschitzian stability, metric regularity, and linear openness; see, e.g., \cite{Mo06a,RoWe98} and the references therein. Recall that a set-valued mapping $\Psi:\R^d\rightrightarrows\R^s$ is {\em Lipschitz-like} around $(\bar z,\bar w)\in\Gr\Psi$ (also known as the pseudo-Lipschitz or Aubin property) with modulus $\kk\ge 0$ if there are neighborhoods $U$ of $\bar z$ and $V$ of $\bar w$ such that
\begin{equation}\label{EqAubinProp}
\Psi(z)\cap V\subset\Psi(u)+\kk\norm{z-u}\B\mbox{ for all }\;z,u\in U.
\end{equation}
The infimum of all such $\kk$ is called the {\em exact Lipschitzian bound} of $\Psi$ around $(\bar z,\bar w)$ and is denoted by $\lip\Psi(\bar z,\bar w)$. If $V=\R^s$, relation \eqref{EqAubinProp} reduces to the (Hausdorff) local Lipschitzian property of set-valued mappings around $\oz$, while in the single-valued case this is nothing else than the classical local Lipschitz continuity of $\Psi$ around the reference point. In terms of \eqref{EqLimCoderiv} we have the robust coderivative characterization of the Lipschitz-like property of $\Psi$ around $(\oz,\ow)$ with the exact bound formula
\begin{eqnarray}\label{cod-cr}
D^*\Psi(\oz,\ow)(0)=\{0\},\quad\lip\Psi(\oz,\ow)=\|D^*\Psi(\oz,\ow)\|
\end{eqnarray}
known as the Mordukhovich criterion \cite[Theorem~9.40]{RoWe98}, where $\|\cdot\|$ stands for the norm of $D^*\Psi(\oz,\ow)$ as a positively homogeneous set-valued mapping from $\R^s$ to $\R^d$. Further, it has been well recognized that $\Psi$ is Lipschitz-like around $(\oz,\ow)$ with modulus $\kk>0$ if and only if its inverse $\Psi^{-1}=M\colon\R^s\tto\R^d$ is {\em metrically regular} around $(\ow,\oz)$ with the same modulus, i.e.,
\begin{eqnarray}\label{EqMetrSubReg}
\dist\big(w;M^{-1}(z)\big)\le\kk\,\dist\big(z;M(w)\big)\;\mbox{ for all }\;z\in U,\;w\in V.
\end{eqnarray}

There are a number of applications for which the robust properties in \eqref{EqAubinProp} and \eqref{EqMetrSubReg} can be relaxed to the weaker ones with putting $u=\oz$ and $z=\oz$ in \eqref{EqAubinProp} and \eqref{EqMetrSubReg}, respectively. The first property is known as {\em calmness} of $\Psi$ at $(\oz,\ow)$, while the second one is known as {\em metric subregularity} of $M$ at $(\ow,\oz)$. Although these properties are equivalent for $\Psi$ and $M=\Psi^{-1}$, we prefer to use metric subregularity in applications to NLPs due to the possibility to formulate it via the initial program data of the original NLP constraint system; see Section~3 for more discussions.

Next we recall two second-order subdifferential constructions for extended-real-valued functions employed below and introduced in the direction initiated in \cite{Mo92}, i.e., by using a coderivative of a first-order subdifferential mapping; this is an appropriate dual set-valued extension of the classical ``derivative-of-derivative" approach to second-order differentiation. Proceeding in this way, we take $f:\R^n\to\oR$, $\xb\in\dom f$, and a basic subgradient $\bar v^\ast\in\partial f(\xb)$ from \eqref{ls} and define the {\em second-order subdifferential} of $f$ at $\xb$ relative to $\bar v^\ast$ as the set-valued mapping $\partial^2f(\xb,\bar v^\ast)\colon\R^n\rightrightarrows\R^n$ given by \cite{Mo92}
\begin{eqnarray}\label{2nd1}
\partial^2f(\xb,\bar v^\ast)(w):=\big(D^\ast\partial f\big)(\xb,\bar v^\ast)(w),\quad w\in\R^n.
\end{eqnarray}
The {\em combined second-order subdifferential} of $f$ at $\ox$ relative to $\ov^*\in\partial f(\ox)$ is defined in this scheme by \cite{MoNg14} replacing the limiting coderivative \eqref{EqLimCoderiv} with its regular counterpart \eqref{EqRegCoderiv}, i.e., by
\begin{eqnarray}\label{2nd2}
\breve{\partial}^2f(\xb,\bar v^\ast)(w):=\big(\widehat D^\ast\partial f\big)(\xb,\bar v^\ast)(w),\quad w\in\R^n.
\end{eqnarray}
When $f$ is $C^2$-smooth around $\xb$, we have $\bar v^\ast=\nabla f(\xb)$ and
$$
\partial^2 f(\xb,\bar v^\ast)(w)=\breve\partial^2 f(\xb,\bar v^\ast)(w)=\big\{\nabla^2f(\ox)w\big\}\;\mbox{ for any }\;w\in\R^n,
$$
where $\nabla^2f(\ox)$ stands for the classical (symmetric) Hessian matrix. Thus the second-order subdifferentials \eqref{2nd1} and \eqref{2nd2} can be treated as the {\em generalized Hessian} constructions.

Now we are ready to formulate and discuss the notion of tilt-stable minimizers for extended-real-valued functions introduced by Poliquin and Rockafellar \cite{PolRo98} without specifying tilt stability moduli and then quantitatively modified and studied by Mordukhovich and Nghia \cite{MoNg14} for the case of given moduli with an explicit calculation of their exact bound.

\begin{Definition} [\bf tilt-stable minimizers]\label{DefTiltStab} Let $f:\R^n\to\oR$, and let $\xb\in\dom f$. Then:

{\bf (i)} $\ox$ is a {\sc tilt-stable local minimizer} of $f$ if there is a number $\gamma>0$ such that the mapping
\begin{equation}\label{EqM_gamma}
M_\gamma(v^\ast):=\argmin\big\{f(x)-\skalp{v^\ast,x}\big|\;x\in\B_\gamma(\xb)\big\},\quad v^*\in\R^n,
\end{equation}
is single-valued and Lipschitz continuous in some neighborhood of $\ov^*=0\in\R^n$ with $M_\gamma(0)=\{\xb\}$.

{\bf (ii)} Given $\kk>0$, the point $\ox$ is a tilt-stable local minimizer of $f$ with {\sc modulus} $\kk$ if there is $\gg>0$ such that $M_\gamma(0)=\{\xb\}$ and the mapping $M_\gg$ in \eqref{EqM_gamma} is single-valued and Lipschitz continuous with modulus $\kk$ around the origin $0\in\R^n$.

{\bf (iii)} The {\sc exact bound of tilt stability} of $f$ at $\ox$ is defined by
\begin{eqnarray}\label{tilt-exact}
\tilt(f,\ox):=\inf_{\gamma>0}\lip M_\gamma(0)
\end{eqnarray}
via the exact Lipschitzian bound of the mapping $M_\gg$ from \eqref{EqM_gamma} around the origin.
\end{Definition}

The main result by Poliquin and Rockafellar \cite[Theorem~1.3]{PolRo98} gives a characterization of tilt-stable minimizers for $f$ in the sense of Definition~\ref{DefTiltStab}(i) via the positive-definiteness of the second-order subdifferential \eqref{2nd1} at the reference point. In this paper we base our investigations on the following quantitative characterization by Mordukhovich and Nghia \cite[Theorem~3.5]{MoNg14}, which uses the combined second-order subdifferential \eqref{2nd2} in a neighborhood of the reference local minimizer and provides, in addition to characterizing tilt-stable minimizers in the sense of Definition~\ref{DefTiltStab}(ii), a precise formula for calculating the exact bound of tilt stability.

\begin{Theorem}[\bf qualitative characterization of tilt-stable minimizers for extended-real-valued function]\label{ThTiltStab}
Let $f:\R^n\to\bar R$ be a l.s.c.\ function having $0\in\partial f(\xb)$. Assume that $f$ is both prox-regular and subdifferentially continuous at $\xb$ for $\bar v^\ast=0$. Then the following assertions are equivalent:

{\bf (i)} The point $\xb$ is a tilt-stable local minimizer of the function $f$ with modulus $\kappa>0$.

{\bf(ii)} There is a constant $\eta>0$ such that
\begin{eqnarray}\label{tilt1}
\skalp{w^\ast,w}\ge\frac 1\kappa\norm{w}^2\;\mbox{ for all }\;w^\ast\in\breve\partial^2f (x,x^\ast)(w),\;(x,x^\ast)\in\Gr\partial f\cap\B_\eta(\xb,0).
\end{eqnarray}
Moreover, the exact bound of tilt stability of $f$ at $\ox$ is calculated by
\begin{eqnarray}\label{tilt2}
\tilt(f,\ox)=\inf_{\eta>0}\sup\left\{\frac{\norm{u}^2}{\skalp{u^\ast,u}}\Big|\;u^\ast\in\breve\partial^2 f (x,x^\ast)(u),\;(x,x^\ast)\in\Gr\partial f\cap\B_\eta(\xb,0)\right\}
\end{eqnarray}
with the convention that $0/0:=0$.
\end{Theorem}

\section{Qualification Conditions in Nonlinear Programming}
\sce

In this section we start a preparatory work for the subsequent second-order characterization of tilt-stable minimizers in NLPs with the system of $C^2$-smooth equality and inequality constraints:
\begin{equation}\label{EqEqualIneqSystem}
\left\{\begin{array}{l}q_i(x)=0\;\mbox{ for }\;i\in E,\\
q_i(x)\le 0\;\mbox{ for }\;i\in I,
\end{array}\right.
\end{equation}
where $E:=\{1,\ldots,l_1\}$ and $I:=\{l_1+1,\ldots,l_1+l_2\}$ are finite index sets for the equality and inequality constraints, respectively.
The main goal of this section is to consider appropriate {\em qualification conditions} needed for characterizing tilt stability in NLPs while being of their own interest.

Denote $l:=l_1+l_2$ and rewrite the constraint system \eqref{EqEqualIneqSystem} in the {\em inclusion form}
\begin{eqnarray}\label{Gamma}
\Gamma:=\big\{x\in\R^n\mv\;q(x)\in\Theta\big\}\;\mbox{ with }\;q(x):=\big(q_1(x),\ldots,q_l(x)\big)\;\mbox{ and }\;\Theta:=\{0\}^{l_1}\times\R^{l_2}_-.
\end{eqnarray}
Consider further the index set of {\em active inequality constraints}
\begin{eqnarray}\label{active}
\I(x):=\big\{i\in I\big|\;q_i(x)=0\big\},\quad x\in\Gamma,
\end{eqnarray}
and for each $x\in\Gamma$ describe the {\em linearized tangent cone} \eqref{EqTan} to $\Gamma$ at this point by
\begin{eqnarray}\label{lin-tan}
\Tlin(x):=\big\{u\in\R^n\big|\;\la\nabla q_i(x),u\ra=0\;\mbox{ for }\;i\in E\;\mbox{ and }\;\la\nabla q_i(x),u\ra\le 0\;\mbox{ for }\;i\in\I(x)\big\}.
\end{eqnarray}
It is easy to observe that the dual/polar cone to \eqref{lin-tan} admits the following representations:
\begin{eqnarray*}
\Tlin(x)^*=\nabla q(x)^*N_\Th\big(q(x)\big)=\disp\Big\{\sum_{\i=1}^l\lm_i\nabla q_i(x)\Big|\;\lambda_i\ge 0\;\mbox{ for }\;i\in\I(x)\;\mbox{ and }\;\lambda_i=0\;\mbox{ for }\;i\in I\setminus\I(x)\Big\}.
\end{eqnarray*}

Recall that the conventional terminology of nonlinear programming understands by ``constraint qualifications" (CQs) any conditions imposed on the constraints of NLPs ensuring that the Lagrange multiplier associated with the cost function in first-order necessary optimality conditions is not zero. For the reader's convenience, let us list the well-recognized CQs, which are compared in what follows with the qualification conditions developed in this paper to study tilt-stable minimizers:\vspace*{0.05in}

$\bullet$  The {\em linear independence constraint qualification} (LICQ) holds at $\xb$ if the gradients of the {\em active constraints} $\{\nabla q_i(\xb)|\;i\in E\cup\I(\xb)\}$ are linearly independent in $\R^n$.

$\bullet$ The {\em Mangasarian Fromovitz constraint qualification} (MFCQ) holds at $\xb$ if if the gradients of the equality constraints $\{\nabla q_i(\xb)|\;i\in E\}$ are linearly independent in $\R^n$ and there exists a vector $u\in\R^n$ satisfying $\la\nabla q_i(\xb),u\ra=0$ for $i\in E$ and $\la\nabla q_i(\xb),u\ra<0$ for $i\in\I(\xb)$.

$\bullet$ The {\em full rank constraint qualification} (FRCQ) holds at $\xb$ if for every subset of the active constraints $\A\subset E\cup\I(\xb)$ we have
\begin{eqnarray*}
\rank\big\{\nabla q_i(\xb)\big|\;i\in\A\big\}=\min\big\{\vert\A\vert,\;n\big\},
\end{eqnarray*}
where $\vert\A\vert$ stands for the cardinality of the set $\A$.

$\bullet$ The {\em constant rank constraint qualification} (CRCQ) holds at $\xb$ if there is a neighborhood $U$ of $\xb$ such that for any index set $\A\subset E\cup\I(\xb)$ the system $\{\nabla q_i(x)|\;i\in\A\}$ has the same rank for all $x\in U$.\vspace*{0.05in}

We have the following implications, which relate the aforementioned CQs at $\ox\in\Gamma$:
\begin{eqnarray}\label{CQ}
LICQ\Longrightarrow MFCQ\;\mbox{ and }\;LICQ\Longrightarrow FRCQ\Longrightarrow CRCQ.
\end{eqnarray}
Indeed, the implications LICQ$\Longrightarrow$MFCQ and LICQ$\Longrightarrow$FRCQ are obvious. The remaining implication FRCQ$\Longrightarrow$CRCQ was observed by Janin (see \cite[Proposition~2.1]{Jan84}) who was the first to introduce and study CRCQ in nonlinear programming.

To proceed further, we recall the equivalent descriptions of MFCQ used in what follows; see, e.g. \cite[Examples~6.40 and 9.44]{RoWe98}. They actually follow from the coderivative characterization \eqref{cod-cr} applied to the inverse of the {\em canonically perturbed} constraint mapping $M_q:\R^n\rightrightarrows\R^l$ defined by
\begin{equation}\label{EqPertMapp}
M_q(x):=q(x)-\Theta,\quad x\in\R^n.
\end{equation}
\begin{Proposition}[\bf equivalent descriptions of MFCQ]\label{ThMFCQ} Given $\xb\in\Gamma$, the validity of MFCQ at $\ox$ is equivalent to each of the following conditions:

{\bf (i)} The mapping $M_q$ is metrically regular around $(\xb,0)$.

{\bf (ii)}$\big[\nabla q(\xb)^\ast\lambda=0,\;\lambda\in N_\Theta\big(q(\xb)\big)\big]\Longrightarrow\lambda=0\in\R^n$.

{\bf (iii)} There exist a positive number $\kappa$ such that
\begin{equation}\label{EqMFCQ_kappa}
\nabla q(\xb)^\ast\lambda\ge\frac{\|\lambda\|}{\kappa}\;\mbox{ for all }\;\lambda\in N_\Theta\big(q(\xb)\big).
\end{equation}
Furthermore, the infimum of the moduli $\kappa$ for which the metric regularity property holds is equal to
\begin{eqnarray*}
\max_{\AT{\lambda\in N_{\Theta}(q(\xb)),}{\norm{\lambda}=1}}\frac 1{\norm{\nabla q(\xb)^\ast\lambda}}.
\end{eqnarray*}
\end{Proposition}

Having in mind the metric regularity description of MFCQ in Proposition~\ref{ThMFCQ}, we define now the following qualification condition, which is clearly weaker than MFCQ and occurs to be very instrumental for the subsequent study of tilt stability.

\begin{Definition}[\bf metric subregularity constraint qualification]\label{DefMetrSubregCQ} Let $\ox\in\Gamma$ for the constraint system \eqref{Gamma}. We say that the {\sc metric subregularity constraint qualification (MSCQ)} holds at $\xb$ if the mapping $M_q$ from \eqref{EqPertMapp} is metrically subregular at $(\xb,0)$.
\end{Definition}

Since in finite-dimensional space all the norms are equivalent, MSCQ can be equivalently described via the existence of a neighborhood $U$ of $\xb$ and a positive number $\kappa$ (for simplicity we keep the same notation for the modulus) such that
\begin{equation}\label{EqMetrSubregRes}
\dist(x;\Gamma)\le\kappa\Big(\sum_{i\in E}
\vert q_i(x)\vert+\sum_{i\in I}\max\big\{q_i(x),0\big\}\Big)\;\mbox{ for all }\;x\in U,
\end{equation}
i.e., for all $x\in U$ sufficiently close to $\xb$ the distance from $x$ to the constraint set $\Gamma$ in \eqref{Gamma} is proportional to the {\em residual} of \eqref{EqEqualIneqSystem} at these points. Note also that the MSCQ property from Definition~\ref{DefMetrSubregCQ} is equivalent to the requirement that the {\em inverse} mapping $S:\R^l\rightrightarrows\R^n$ given by
\begin{eqnarray}\label{S}
S(y):=M_q^{-1}(y)=\big\{x\in\R^n\big|\;y\in q(x)-\Theta\big\},\quad y\in\R^l,
\end{eqnarray}
is {\em calm} at $(0,\xb)$. We prefer to deal primarily with MSCQ instead of the calmness requirement on $S$ due to the fact that the condition in \eqref{EqMetrSubregRes} is formulated in terms of the {\em initial program data} $q_i$ while the inverse mapping $S$ may not be in hands, and it is usually hard to construct it.

Observe that the imposed MSCQ is indeed a {\em constraint qualification} at $\ox$ in the standard sense of nonlinear programming recalled above. Indeed, the following implication follows from \cite[Proposition~1]{HenOut05} and the aforementioned relationship between metric subregularity and calmness:
\begin{eqnarray}\label{aba}
\big[MSCQ\;\mbox{ at }\;\ox\in\Gamma\big]\Longrightarrow\big[T_\Gamma(\xb)=\Tlin(\xb)\big],
\end{eqnarray}
where the right-hand side equality in \eqref{aba}, saying that the tangent cone \eqref{EqTan} to $\Gamma$ at $\ox$ agrees with the one \eqref{lin-tan} to the linearized constraints, is known as the {\em Abadie constraint qualification} (ACQ) for \eqref{EqEqualIneqSystem} at $\ox$, which a CQ in the standard NLP sense.\vspace*{0.05in}

In order to conduct our subsequent analysis of tilt stability for local minimizers in NLPs, MSCQ alone is not enough. As the reader can see below, just one additional qualification condition on the constraint system \eqref{EqEqualIneqSystem} at $\ox$ is needed. To define this new condition, let us first introduce some objects associated with \eqref{EqEqualIneqSystem} and \eqref{Gamma}. Given vectors $x\in\Gamma$ and $x^\ast\in\R^n$, consider the {\em set of multipliers}
\begin{eqnarray}\label{Lambda}
\Lambda(x,x^\ast):=\big\{\lambda\in N_\Theta\big(q(x)\big)\big|\;\nabla q(x)^\ast\lambda=x^\ast\big\}
\end{eqnarray}
and the corresponding collection of {\em strict complementarity indexes}
\begin{eqnarray}\label{I+}
I^+(\lambda):=\big\{i\in I\big|\;\lambda_i>0\big\}\;\mbox{ for }\;\lm=(\lm_1,\ldots,\lm_l)\in\Th^*.
\end{eqnarray}
Denote by $\E(x,x^\ast)$ the collection of all the {\em extreme points} of the closed and convex set of multipliers $\Lambda(x,x^\ast)$ and recall that $\lambda\in\Lambda(x,x^*)$ belongs to $\E(x,x^\ast)$ if and only if the family of gradients $\{\nabla q_i(x)\mv i\in E\cup I^+(\lambda)\}$ is linearly independent. Hence $\E(x,x^*)\ne\emp$ if and only if $\Lambda(x,x^*)\ne\emp$ and the gradients of the equality constraints $\{\nabla q_i(x)|\;i\in E\}$ are linearly independent.\vspace*{0.05in}

Now we are ready to introduce the new qualification condition for the constraint system \eqref{EqEqualIneqSystem}.

\begin{Definition}[\bf bounded extreme point property]\label{DefBEPP} We say that the {\sc bounded extreme point property (BEPP)} holds at $\xb$ if the gradients of the equality constraints $\{\nabla q_i(\xb)|\;i\in E\}$ are linearly independent and there exist a neighborhood $U$ of $\xb$ and a number $\kappa>0$ such that
\begin{equation}\label{EqBEPP}
\E(x,x^\ast)\subset\kappa\norm{x^\ast}\B\;\mbox{ for all }\;x\in\Gamma\cap U\;\mbox{ and }\;x^\ast\in\R^n.
\end{equation}
\end{Definition}

In contrast to the case of MSCQ, we do not claim that BEPP is a constraint qualification in the standard sense. Therefore the term ``{\em qualification condition}" seems to be appropriate for both MSCQ and BEPP. In what follows we study the properties MSCQ and BEPP simultaneously and apply them together to deriving pointbased conditions for tilt stable minimizers in NLPs.

The next proposition shows, in particular, that {\em each} of the constraint qualifications MFCQ and CRCQ, and thus the stronger ones from \eqref{CQ}, ensures the validity of {\em both} MSCQ and BEPP.

\begin{Proposition}[\bf robustness of MSCQ and BEPP with sufficient conditions for their validity]\label{LemSuffCondBEPP} Given $\ox\in\Gamma$ from \eqref{Gamma}, the following assertions hold:

{\bf (i)}  Both MSCQ and BEPP are {\sc robust} properties in the sense that if either MSCQ or BEPP is satisfied at the reference point $\ox\in\Gamma$ that there is a neighborhood $U$ of $\ox$ such that the corresponding property is satisfied at any point $x\in\Gamma\cap U$.

{\bf (ii)} Let either MFCQ or CRCQ hold at $\xb\in\Gamma$ and that $($in the case of CRCQ$)$ the gradients of the equality constraints $\{\nabla q_i(\xb)|\;i\in E\}$ are linearly independent in $\R^n$. Then both the qualification conditions MSCQ and BEPP are satisfied at $\xb$.
\end{Proposition}
{\bf Proof.} Assertion (i) for both MSCQ and BEPP follows directly from the definitions of these qualification conditions. Also it is straightforward to deduce from the the characterization of MFCQ in Proposition~\ref{ThMFCQ}(i) that the validity of MFCQ at $\ox$ implies that MSCQ holds at this point.

Suppose now that CRCQ holds at $\xb$. Then it follows from \cite[Proposition~2.5]{Jan84} that the mapping $S\colon\R^l\tto\R^n$ from \eqref{S}
is calm at $(0,\xb)$, and hence $M_q(x)=S^{-1}(x)$ is metrically subregular at $(\xb,0)$.

Further, let us check that MFCQ at $\ox$ yields the validity of BEPP at this point. Using the equivalent description of MFCQ in Proposition~\ref{ThMFCQ}(iii), take $\kk>0$ from \eqref{EqMFCQ_kappa} and find a neighborhood $U$ of $\xb$ such that $\I(x)\subset\I(\xb)$ and that the estimate $\norm{\nabla q(x)-\nabla q(\xb)}\le\frac 1{2\kappa}$ holds on $\Gamma\cap U$. Then for every $x\in\Gamma\cap U$, $x^\ast\in\R^n$, and $\lambda\in\E(x,x^\ast)$ we get $\lambda\in N_\Theta(q(x))\subset N_\Theta(q(\xb))$ and hence
\begin{eqnarray*}
\norm{x^\ast}=\norm{\nabla q(x)^\ast\lambda}\ge\norm{\nabla q(x)^\ast\lambda}-\norm{\nabla q(x)-\nabla q(\xb)}\cdot\norm{\lambda}\ge\frac{\norm{\lambda}}{\kappa}-\frac{\norm{\lambda}}{2\kappa}=\frac{\norm{\lambda}}{2\kappa}.
\end{eqnarray*}
This shows that $\E(x,x^\ast)\subset 2\kappa\norm{x^\ast}\B$ and thus justifies that BEPP holds at $\ox$.

It remains to verify that the validity of CRCQ at $\xb$ together with the linear independence of $\{\nabla q_i(\xb)|\;i\in E\}$ implies that BEPP holds at this point. Assuming the contrary and employing the imposed linear independence allow us to find sequences $x_k\setto\Gamma\xb$, $x_k^\ast\in\R^n$, and $\lambda^k\in\E(x_k,x_k^\ast)$ such that
$$
\nabla q(x_k)^\ast\frac{\lambda^k}{\|\lambda^k\|}\to 0\;\mbox{ as }\;k\to\infty.
$$
Passing to a subsequence if necessary gives us a vector $\lm\in\R^n$ with $\|\lm\|=1$ and such that $\lambda^k/\|\lm^k\|\to\lambda$ as $k\to\infty$. Since $\nabla q(\xb)^\ast\lambda=0$, it follows from the above that the gradient family $\{\nabla q_i(\xb)|\;i\in E\cup I^+(\lm)\}$ with $I^+(\lm)$ from \eqref{I+} is {\em linearly dependent}, i.e., consists of linearly dependent vectors in $\R^n$. Then for each $i\in I^+(\lm)$ we have $\lambda^k_i>0$ whenever $k\in\N$ is sufficiently large, which shows that $i\in\I(x_k)\subset\I(\xb)$ according to \eqref{active} and justifies in turn that
$I^+(\lm)\subset I^+(\lm^k)\subset\I(\xb)$. Now the assumed CRCQ at $\ox$ ensures that the family $\{\nabla q_i(x_k)|\;i\in E\cup I^+(\lm)\}$ is linearly dependent, and hence the family $\{\nabla q_i(x_k)|\;i\in E\cup I^+(\lambda^k)\}$ is linearly dependent as well. This clearly contradicts, due to the discussion right before Definition~\ref{DefBEPP}, that $\lambda^k\in\E(x_k,x_k^\ast)$ for large numbers $k\in\N $. Thus BEPP holds at $\xb$, which completes the proof of the proposition. $\h$

\section{MSCQ and BEPP via Second-Order Analysis}
\sce

In this section we employ second-order derivatives of the constraint functions $q_i$ at $\ox$ to effectively support MSCQ and BEPP and also use these qualification conditions to describe some second-order constructions of variational analysis needed in what follows.

Note that LICQ, MFCQ, and FRCQ are {\em pointbased} conditions in contrast to CRCQ, MSCQ, and BEPP that involve {\em neighborhood} points in their definitions.  It is worth mentioning to this end that the papers by Gfrerer \cite{Gfr11} and by Li and Mordukhovich \cite{LM12} contain some (different) pointbased sufficient conditions for metric subregularity of general set-valued mappings, which are ``almost necessary" for this property. On the other hand, in the same paper \cite{Gfr11} Gfrerer introduced the pointbased  {\em second-order sufficient condition for metric subregularity} (SOSCMS) formulated below in Theorem~\ref{ThSSOSCMS} that allowed him to derive ``no-gap" second-order necessary and sufficient conditions for metric subregularity, where the difference between the necessity and sufficiency is the change from inequality to strict inequality. The reader can find further applications of this condition in the recent papers \cite{Gfr14a,GfrKl14,GfrOut14a,GfrOut14b}.

The next theorem shows that SOSCMS ensures the validity of not only MSCQ but also of BEPP, i.e., of {\em both} qualification conditions we use for our subsequent pointbased characterizations of tilt stability in NLPs. In fact, this theorem provides a stronger version of the aforementioned results. Namely, we consider the case when the constraint system \eqref{EqEqualIneqSystem} can be split into the following two subsystems with both equality and inequality constraints:
\begin{eqnarray*}
\begin{array}{ll}
q_i(x)=0\;\mbox{ for }\;i\in E_1:=\big\{1,\ldots,l_1^1\big\}\;\mbox{ and }\;q_i(x)=0\;\mbox{ for }\;i\in E_2:=\big\{l_1^1+1,\ldots,l_1^1+l_1^2=l_1\},\\\\
q_i(x)\le0\;\mbox{ for }\;i\in I_1:=\big\{1,\ldots,l_2^1\big\}\;\mbox{ and }\;q_i(x)\le 0\;\mbox{ for }\;i\in I_2:=\big\{l_2^1+1,\ldots,l_2^1+l_2^2=l_2\big\}
\end{array}
\end{eqnarray*}
in such a way that it is known in advance that for the second system $q_i(x)=0$ as $i\in E_2$ and $q_i(x)\le 0$ as $i\in I_2$ both MSCQ and BEPP are surely satisfied. In particular, it happens by Proposition~\ref{LemSuffCondBEPP}(ii) if CRCQ is fulfilled at $\ox$ and the gradient vectors $\{\nabla q_i(\xb)|\;i\in E_2\}$ are linearly independent. One of the reasons for this is that, although in the absence of FRCQ no pointbased conditions for verifying CRCQ is known in terms of the gradients $\nabla g_i(\ox)$, there exist other easily verifiable conditions that ensure the validity of CRCQ without using any derivatives. For instance, it is well known that CRCQ holds at every $\ox\in\Gamma$ if the functions $q_i$ are {\em linear} whenever $i\in E_2\cup I_2$. Note that if we ignore the $(E_2,I_2)$-system in the following theorem (i.e., put $E_2=I_2=\emp$), then it merely asserts the validity of both MSCQ and BEPP for the original constraint system \eqref{EqEqualIneqSystem} at $\ox$ under the pointbased SOSCMS assumption imposed at this point, which is surely implied by MFCQ.

\begin{Theorem}[\bf MSCQ and BEPP from SOSCMS]\label{ThSSOSCMS} Let the gradients of the equality constraints $\{\nabla q_i(\xb)|\;i\in E\}$ be linearly independent, and the system $q_i(x)=0$ for $i\in E_2$ and $q_i(x)\le 0$ for $i\in I_2$ fulfill both MSCQ and BEPP at $\xb\in\Gamma$. Impose further the following SOSCMS at $\ox$: for every vector $0\not=u\in\Tlin(\xb)$ from the linearized constraint tangent cone \eqref{lin-tan} we have the implication
\begin{eqnarray}\label{SOSCMS}
\big[\lambda\in N_\Theta\big(q(\xb)\big),\;\nabla q(\xb)^\ast\lambda=0,\;\big\la u,\nabla^2\la\lambda,q\ra(\xb)u\big\ra\ge
0\big]\Longrightarrow\disp\sum_{i\in E_1\cup I_1}\vert\lambda_i\vert=0.
\end{eqnarray}
Then both MSCQ and BEPP are satisfied for the original constraint system \eqref{EqEqualIneqSystem} at $\xb$.
\end{Theorem}
{\bf Proof.} Observe first that the implication
$$
SOSCMS\Longrightarrow MSCQ\;\mbox{ at }\;\ox\in\Gamma
$$
in the general setting of Theorem~\ref{ThSSOSCMS} follows from the combination of Theorem~2.6 and Lemma~2.7 in \cite{Gfr14a}. It remains to verify the other implication
\begin{eqnarray}\label{SOS-BEPP}
SOSCMS\Longrightarrow BEPP\;\mbox{ at }\;\ox\in\Gamma.
\end{eqnarray}

Assuming the contrary to \eqref{SOS-BEPP}, find sequences $x_k\setto{\Gamma}\xb$ and $\lambda^k\in N_\Theta(q(x_k))$ so that the gradients $\{\nabla q_i(x_k)|\;i\in E\cup I^+(\lambda^k)\}$ are linearly independent in $\R^n$ and that
$$
\|\nabla q(x_k)^\ast\lambda^k\|<\norm{\lambda^k}/k\;\mbox{ for all }\;k\in\N.
$$
Passing to a subsequence of $k\to\infty$ gives us $\lambda\in N_\Theta(q(\xb))$ with $\|\lm\|=1$ such that $\lambda^k/\norm{\lambda^k}\to\lambda$
and $\nabla q(\xb)^\ast\lambda=0$, which yields the linear dependence of the gradients $\{\nabla q_i(\xb)\mv i\in E\cup I^+(\lambda)\}$.
Since $I^+(\lambda)\subset I^+(\lambda^k)$ for the index sets \eqref{I+} and the family of gradients $\{\nabla q_i(x_k)|\;i\in E\cup I^+(\lambda^k)\}$ is linearly independent, it follows that $x_k\not=\xb$ for all $k$ sufficiently large. Passing to a subsequence again allows us to find $u\in\R^n$ with $\|u\|=1$ for which $(x_k-\xb)/\norm{x_k-\xb}\to u$ as $k\to\infty$. We obviously have for any active constraint $i\in E\cup\I(\xb)$ that
\[
\la\nabla q_i(\xb),u\ra=\lim_{k\to\infty}\frac{q_i(x_k)-q_i(\xb)}{\norm{x_k-\xb}}\begin{cases}=0&\mbox{ if }\;i\in E,\\\le 0&\mbox{ if }\;i\in\I(\xb)
\end{cases}
\]
showing that $\nabla q(\xb)u\in T_\Theta(q(\xb))$ and consequently that $u\in\Tlin(\xb)$. Furthermore, we deduce from
$I^+(\lambda)\subset I^+(\lambda^k)$ for large $k$ that the condition $\la\lambda,q(x_k)\ra=0$ yields
\begin{eqnarray*}
0=\lim_{k\to\infty}\frac{\la\lambda,q(x_k)-q(\xb)\ra}{\norm{x_k-\xb}^2}=\frac 12\big\la u,\nabla^2\la\lambda,q\ra(\xb)u\big\ra.
\end{eqnarray*}
Let us check now that $\sum_{i\in E_1\cup I_1}\vert\lambda_i\vert>0$, which clearly contradicts the SOSCMS assumption in \eqref{SOSCMS}.
Indeed, setting $\lambda^{k,j}:=(\lambda^k_i)_{i\in E_j\cup I_j}$ and $q^{(j)}:=(q_i)_{i\in E_j\cup I_j}$ for $j=1,2$, we get
\begin{eqnarray*}
\frac{\|\lambda^{k,2}\|}{\kappa_2}-\big\|\nabla q^{(1)}(x_k)^\ast\lambda^{k,1}\big\|
&\le&\big\|\nabla q^{(2)}(x_k)^\ast\lambda^{k,2}\big\|-\big\|\nabla q^{(1)}(x_k)^\ast\lambda^{k,1}\big\|\\
&\le&\big\|\nabla q^{(2)}(x_k)^\ast\lambda^{k,2}+\nabla q^{(1)}(x_k)^\ast\lambda^{k,1}\big\|\\
&=&\big\|\nabla q(x_k)^\ast\lambda^k\big\|<\frac{\|\lambda^k\|}{k},
\end{eqnarray*}
where $\kappa_2$ denotes a positive number such that the BEPP assumption \eqref{EqBEPP} holds for the ``second" constraint subsystem $q_i(x)=0$ as $i\in E_2$ and $q_i(x)\le 0$ as $i\in I_2$. Dividing the latter inequality by $\norm{\lambda^k}$ and passing to the limit as $k\to\infty$, we obtain
\begin{eqnarray*}
\big\|\nabla q^{(1)}(\xb)^\ast\lambda^{(1)}\big\|\ge\frac{\norm{\lambda^{(2)}}}{\kappa_2},\;\mbox{ where }\;\lambda^{(j)}:=(\lambda_i)_{i\in E_j\cup I_j}\;\mbox{ for }\;j=1,2.
\end{eqnarray*}
Since $\norm{\lambda}=\norm{(\lambda^{(1)},\lambda^{(2)})}=1$ we easily conclude that $\lambda^{(1)}\not=0$, and thus $\sum_{i\in E_1\cup I_1}\vert\lambda_i\vert>0$. The obtained contradiction with \eqref{SOSCMS} justifies \eqref{SOS-BEPP} and completes the proof of the theorem. $\h$\vspace*{0.05in}

Next we show that the simultaneous validity of MSCQ and BEPP at $\ox$ ensures that the indicator function $\dd_\Gamma$ of the constraint set $\Gamma$ from \eqref{Gamma} belongs to the basic in second-order analysis class of prox-regular and subdifferentially continuous functions needed for the second-order characterization of tilt-stable minimizers in the abstract extended-real-valued setting of Theorem~\ref{ThTiltStab}.

\begin{Proposition}[\bf prox-regularity and subdifferential continuity of the constraint indicator function]\label{LemProxRegDelta} Assume that both MSCQ and BEPP hold at $\xb\in\Gamma$. Then there is a neighborhood $U$ of $\xb$ such that for all $x\in\Gamma\cap U$ we have the equalities
\begin{eqnarray}\label{prox1}
\partial\delta_\Gamma(x)=\Hat\partial\delta_\Gamma(x)=\nabla q(x)^*N_\Theta\big(q(x)\big).
\end{eqnarray}
Furthermore, $\delta_\Gamma$ is prox-regular and subdifferentially continuous at $\xb$ for every $\xba\in\partial\delta_\Gamma(\xb)$.
\end{Proposition}
{\bf Proof.} The validity of the second equality in \eqref{prox1} for $x=\ox$ under MSCQ at this point follows from implication \eqref{aba} and the fact that ACQ at $\ox$ immediately implies the dual condition $\Hat N_\Gamma(\ox)=\Tlin(\ox)^*(=\nabla q(x)^*N_\Theta(q(\ox)))$ known as the {\em Guignard constraint qualification} (GCQ). Hence Proposition~\ref{LemSuffCondBEPP}(i) ensures the second equality in \eqref{prox1} for all $x\in\Gamma$ near $\ox$. We now show that the first equality in \eqref{prox1} is also satisfied if in addition BEPP holds at $\ox$ and hence around this point.

Suppose by Definition~\ref{DefBEPP} that \eqref{EqBEPP} holds with some $\kappa$ and that for all $x\in U$ the gradients $(\nabla q_i(x))_{i\in E}$ are linearly independent. Fix $x\in\Gamma\cap U$ and $x^\ast\in\partial\delta_\Gamma(x)$ and then find sequences $x_k\setto{\Gamma} x$ and $x_k^\ast\to x^\ast$ with $x_k^\ast\in\hat\partial\delta_\Gamma(x_k)$ for all $k$. Since $x_k\in U$ when $k$ is large, we get $\Lambda(x_k,x_k^\ast)\ne\emp$ for the set of multipliers and consequently $\E(x_k,x_k^\ast)\ne\emp$ for the collection of its extreme points. Picking $\lambda^k\in\E(x_k,x_k^\ast)$ for each $k$ gives us $\norm{\lambda^k}\le\kappa\norm{x_k^\ast}$. Thus the sequence $\{\lambda_k\}$ is bounded and converges therefore to some $\lambda\in\R^l$ along a subsequence. We obviously have $\lambda\in N_\Theta(q(x))$ and $\nabla q(x)^\ast\lambda=x^\ast$ showing that $x^\ast\in\hat\partial\delta_\Gamma(x)$. Since the opposite inclusion $\hat\partial\delta_\Gamma(x)\subset\partial\delta_\Gamma(x)$ always holds, it tells us that $\partial\delta_\Gamma(x)=\hat\partial\delta_\Gamma(x)$ for every $x\in\Gamma\cap U$ proving in this way the first equality in \eqref{prox1}.

Considering the last statement of the proposition, observe easily from the definitions that $\delta_\Gamma$ is subdifferentially continuous at $\xb$ for $\xba\in\partial\delta_\Gamma(\xb)$ and that the epigraph of $\delta_\Gamma$ is closed, i.e., $\delta_\Gamma$ is l.s.c.\ on $\R^n$. Taking now $u,x\in\Gamma\cap U$ and $x^\ast\in\partial\delta_\Gamma(x)$, pick $\lambda\in\E(x,x^\ast)$ and get by BEPP \eqref{EqBEPP} that $\norm{\lambda}\le\kappa\norm{x^\ast}$. Since $\la\lambda,q(x)\ra=0$, $\nabla q(x)^\ast\lambda=x^\ast$, and $\lambda_i q_i(u)\le 0$ for $i\in E\cup I$, we conclude that
\begin{eqnarray*}
\delta_\Gamma(u)-\delta_\Gamma(x)&=0&\ge\la\lambda,q(u)-q(x)\ra\ge\la\lambda,\nabla q(x)(u-x)\ra-\frac{\gg}2\norm{\lambda}\cdot\norm{u-x}^2
\\&\ge&\skalp{x^\ast,u-x}-\frac{\gg}2\kappa\norm{x^\ast}\cdot\norm{u-x}^2,
\end{eqnarray*}
where $\gg:=\sup\big\{\norm{\nabla^2q(y)}\;\big|\;y\in U\big\}$. This verifies the prox-regularity \eqref{prox-reg} of $\dd_\Gamma$ at $\xb$ for every subgradient $\xba\in\partial\delta_\Gamma(\xb)$ and thus completes the proof of the proposition. $\h$\vspace*{0.05in}

As indicated by one of the referees, the prox-regularity of $\delta_\Gamma$ under MSCQ in Proposition~\ref{LemProxRegDelta} can be derived from \cite[Theorem~31(b)]{CT10}, although the notion of MSCQ was not defined therein.\vspace*{0.03in}

To proceed further, recall the definition of the {\em critical cone} to $\Gamma$ at $(x,x^\ast)\in\Gr\hat\partial\delta_\Gamma$ given by
\begin{eqnarray}\label{crit-cone}
K(x,x^\ast):=T_\Gamma(x)\cap\{x^\ast\}^\perp
\end{eqnarray}
via the tangent cone \eqref{EqTan} and define the {\em multiplier set in a direction} $v\in\R^n$ by
\begin{eqnarray}\label{crit-mult}
\Lambda(x,x^\ast;v):=\argmax\big\{\big\la v,\nabla^2\la\lambda,q\ra(x)v\big\ra\big|\;\lambda\in\Lambda(x,x^\ast)\big\}.
\end{eqnarray}
Note that \eqref{crit-mult} consists of {\em optimal solutions} to a {\em linear program} over the feasible set of multipliers \eqref{Lambda}. This ``critical multiplier set" plays a crucial role in our subsequent study of tilt stability.\vspace*{0.05in}

The next proposition collects some properties of sets \eqref{crit-cone} and \eqref{crit-mult} needed in what follows.

\begin{Proposition}[\bf critical cone and multipliers in critical directions under MSCQ]\label{LemElemProp} Let MSCQ hold at $\ox\in\Gamma$, and let $(x,x^\ast)\in\Gr\hat\partial\delta_\Gamma$ be any pair such that $x\in\Gamma$ is sufficiently close to $\ox$. Then the following assertions are satisfied:

{\bf(i)} For every multiplier $\lambda\in\Lambda(x,x^\ast)$ we have
\[
K(x,x^\ast)=\left\{v\in\R^n\Big|\;\la\nabla q_i(x),v\ra\begin{cases}=0&\mbox{if }\;i\in E\cup I^+(\lambda)\\\le 0&\mbox{if }\;i\in
\I(x)\setminus I^+(\lambda)\end{cases}\right\}.
\]

{\bf(ii)} There exist a multiplier $\tilde\lambda\in\Lambda(x,x^\ast)$ with $I^+(\tilde\lambda)=I^+:=\bigcup_{\lambda\in\Lambda(x,x^\ast)}
I^+(\lambda)$ and some vector $v\in K(x,x^\ast)$ satisfying the conditions
\[
\la\nabla q_i(x),v\ra\begin{cases}=0&\mbox{if }\;i\in E\cup I^+(\Tilde\lambda),\\<0&\mbox{if }\;i\in
\I(x)\setminus I^+(\Tilde\lambda).
\end{cases}
\]

{\bf(iii)} For every vector $v\in K(x,x^\ast)$ we have $\Lambda(x,x^\ast;v)\ne\emp$.
\end{Proposition}
{\bf Proof.} To verify the first assertion, we use the robustness of MSCQ by Proposition~\ref{LemSuffCondBEPP}(i) and implication \eqref{aba}, which ensure that $\Tlin(x)=T_\Gamma(x)$ for all $x\in\Gamma$ around $\ox$. This yields $\Lambda(x,x^\ast)\ne\emp$ and thus (i) follows from the observation that $v\in K(x,x^\ast)$ if and only if
$$
v\in\Tlin(x)\;\mbox{ and }\;0=\skalp{x^\ast,v}=\sum_{i\in E\cup I}\lambda_i\la\nabla q_i(x),v\ra=\sum_{i\in E\cup I^+(\lambda)}\lambda_i\la\nabla q_i(x),v\ra,
$$
where the multiplier $\lambda\in\Lambda(x,x^\ast)$ is chosen arbitrarily. Assertion (ii) follows from \cite[Lemma~2]{GfrOut14a}. To justify finally assertion (iii), we employ the {\em dual} second-order necessary condition for metric subregularity from \cite[Theorem~6.1]{Gfr11} and obtain in this way that for every $v\in\Tlin(x)\supset K(x,x^\ast)$ and every $\lambda\in N_\Theta(q(x))$ with $
\nabla q(x)^\ast\lambda=0$ it follows that $\la v,\nabla^2\la\lambda,q\ra(x)v\ra\le 0$. The latter inequality implies that the linear optimization problem
\begin{eqnarray}\label{lp}
\mbox{maximize }\;\big\la v,\nabla^2\la\lambda,q\ra(x)v\big\ra\;\mbox{ subject to }\;\lambda\in\Lambda(x,x^*)
\end{eqnarray}
admits an optimal solution, which exactly means that $\Lambda(x,x^\ast;v)\ne\emp$. $\h$

\section{Calculations of Second-Order Generalized Derivatives for NLPs}
\sce

In this section we present {\em precise calculations} of some generalized second-order derivative constructions for the indicator function $\dd_\Gamma$ of the feasible solution set \eqref{Gamma} given by equality and inequality constraints via the second-order derivatives of the constraint functions $q_i$ as well as the critical cone \eqref{crit-cone} and the set of multipliers in critical directions \eqref{crit-mult}. The theorem below extends the recent results by Gfrerer and Outrata \cite{GfrOut14a} regarding the following major issues:\vspace*{0.05in}

$\bullet$ It concerns not only inequality but also equality constraints in \eqref{EqEqualIneqSystem}.

$\bullet$ It replaces a certain relaxed uniform metric regularity property in the vicinity of the reference point employed in \cite{GfrOut14a} by the weaker BEPP qualification condition imposed at this point.

\begin{Theorem}[\bf generalized second-order derivatives of the constraint indicator function under MSCQ and BEPP]\label{ThTnaConeGrNormalCone}
Given $\ox\in\Gamma$, assume that both MSCQ and BEPP hold at $\ox$. Then for any $x\in\Gamma$ sufficiently close to $\ox$ and any regular subgradient $x^*\in\Hat\partial\dd_\Gamma(x)$ the following assertions hold:

{\bf (i)} The tangent cone \eqref{EqTan} to the graph of $\Hat\partial\dd_\Gamma$ is calculated by
\begin{equation}\label{EqTanConeGrNormalCone}
T_{{\rm gph}\,\Hat\partial\delta_\Gamma}(x,x^\ast)=\big\{(v,v^\ast)\in\R^{2n}\big|\;\exists\,\lambda\in\Lambda(x,x^\ast;v)\;\mbox{ with }\;
v^\ast\in\nabla^2\la\lambda,q\ra(x)v+N_{K(x,x^\ast)}(v)\big\}.
\end{equation}

{\bf (ii)} Assume that $\Lambda(x,x^\ast;\cdot)$ is constant on $K(x,x^\ast)\setminus\{0\}$ and take an arbitrary multiplier
\begin{eqnarray}\label{mult1}
\lm\in\left\{\begin{array}{ll}\Lambda(x,x^\ast;v)\;\mbox{ for }\;v\in K(x,x^*)\setminus\{0\}&\mbox{ if }\;K(x,x^*)\ne\{0\},\\
\Lambda(x,x^*;0)&\mbox{ otherwise}.
\end{array}\right.
\end{eqnarray}
Then we have the simplified tangent cone formula
\begin{equation}\label{EqTangConeConstLambda}
T_{{\rm gph}\,\Hat\partial\delta_\Gamma}(x,x^\ast)=\big\{(v,v^\ast)\in\R^{2n}\big|\;v^\ast\in\nabla^2\la\lambda,q\ra(x)v+N_{K(x,x^\ast)}(v)\big\}.
\end{equation}
Furthermore, the regular normal cone \eqref{fn1} to the graph of $\Hat\partial\dd_\Gamma$ is calculated by
\begin{equation}\label{EqNormalConeConstLambda}
\widehat N_{{\rm gph}\,\Hat\partial\delta_\Gamma}(x,x^\ast)=\big\{(w^\ast,w)\in\R^{2n}\big|\;
w\in K(x,x^\ast),\;w^\ast\in-\nabla^2\la\lambda,q\ra(x)w+K(x,x^\ast)^*\big\}
\end{equation}
with an arbitrary multiplier $\lm$ taken from \eqref{mult1}.
\end{Theorem}
{\bf Proof.} Consider the equivalent representation of $\Gamma$ obtained by replacing the equality constraints with two inequalities, i.e., the following one:
\begin{equation}\label{EqIneqSystem}
\tilde q_i(x)\le 0\;\mbox{ for }\;i=1,\ldots,\tilde l,\;\mbox{ where }\;\tilde l:=2l_1+l_2,
\end{equation}
\[\tilde q_i(x):=\begin{cases}q_i(x)&\mbox{if $1\le i\le l_1$,}\\
-q_{i-l_1}(x)&\mbox{if $l_1+1\le i\le 2l_1$,}\\
q_{i-l_1}(x)&\mbox{if $2l_1+1\le i\leq 2l_1+l_2$.}\end{cases}\]
It is easy to conclude from the metric subregularity description \eqref{EqMetrSubregRes} that the modified constraint mapping $M_{\tilde q}(x):=\tilde q(x)-\R^{\tilde l}_-$ for \eqref{EqIneqSystem} is metrically subregular at $(\ox,0)$ if and only the original one $M_q$ from \eqref{EqPertMapp} has the same property at this point. Proceeding similarly to the proof of \cite[Theorem~1]{GfrOut14a} while using the BEPP condition at $\ox$ (and hence at points nearby), we arrive at the representation
\begin{equation}
\label{EqTanConeGrNormalConeInequ}
T_{{\rm gph}\,\hat\partial\delta_\Gamma}(x,x^\ast)=\big\{(v,v^\ast)\in\R^{2n}\big|\;\exists\,\tilde\lambda\in\tilde\Lambda(x,x^\ast;v)\;\mbox{ with }\;v^\ast\in\nabla^2\la\tilde\lambda,q\ra(x)v+N_{K(x,x^\ast)}(v)\big\},
\end{equation}
where the sets of multipliers $\tilde\Lambda(x,x^\ast)$ and $\tilde\Lambda(x,x^\ast;v)$ are defined as in \eqref{Lambda} and \eqref{crit-mult}, respectively, but for the extended inequality system \eqref{EqIneqSystem}, i.e.,
$$
\tilde\Lambda(x,x^\ast):=\big\{\tilde\lambda\in N_{R^{\tilde l}_-}\big(\tilde q(x)\big)\big|\;\nabla\tilde q(x)^\ast\tilde\lambda=x^\ast\big\},\;\tilde\Lambda(x,x^\ast;v):=\argmax\big\{\big\la v,\nabla^2\la\tilde\lambda,\tilde q\ra(x)v\big\ra\big|\;\tilde\lambda\in\tilde\Lambda(x,x^\ast)\big\}.
$$
The only essential difference from the proof of \cite[Theorem~1]{GfrOut14a} is that now we need to justify the following fact under the MSCQ and BEPP assumptions made: for any sequences $x_k\setto\Gamma x$ and $x_k^\ast\to x^\ast$ with $x_k^\ast\in\Hat N_\Gamma(x_k)$ there exists a {\em bounded} sequence of multipliers $\tilde\lambda^k\in\tilde\Lambda(x_k,x_k^\ast)$ as $k\in\N$. To verify this, observe that BEPP together with MSCQ guarantees the existence of a sequence $\lambda^k\in\E(x_k,x_k^\ast)$ satisfying the estimate $\norm{\lambda^k}\le\kappa\norm{x_k^\ast}$ with some constant $\kappa>0$ independent of $k\in\N$. Defining
\begin{eqnarray*}
\begin{array}{ll}
&\tilde\lambda_{i+l_1}^k:=\lambda_i^k\;\mbox{ for }\;i=l_1+1,\ldots,l,\\
&\tilde\lambda_i^k:=\lambda_i^k,\;\tilde\lambda_{i+l_1}^k=0\;\mbox{ for }\;i=1,\ldots,l_1\;\mbox{ with }\;\lambda_i^k\ge 0,\\
&\tilde\lambda_i^k:=0,\;\tilde\lambda_{i+l_1}^k:=-\lambda_i^k\;\mbox{ for }\;i=1,\ldots,l_1\;\mbox{ with }\;\lambda_i^k<0,
\end{array}
\end{eqnarray*}
we get $\tilde\lambda^k\in\tilde\Lambda(x_k,x^\ast_k)$ and $\norm{\tilde\lambda^k}=\norm{\lambda^k}$. It shows that the sequence $\{\tilde\lambda^k\}$ is bounded, and so formula \eqref{EqTanConeGrNormalConeInequ} holds by \cite[Theorem~1]{GfrOut14a}. It is easy to see that the set on the right-hand side of \eqref{EqTanConeGrNormalCone} is the same as the one on the right-hand side of \eqref{EqTanConeGrNormalConeInequ}, which thus verifies the claimed representation \eqref{EqTanConeGrNormalCone}.

To justify \eqref{EqTangConeConstLambda} and its dual version \eqref{EqNormalConeConstLambda}, note that for every pair $0\not=v_1,v_2\in K(x,x^\ast)$ we have $\Lambda(x,x^\ast;v_1)=\Lambda(x,x^\ast;v_2)$ if and only if $\tilde\Lambda(x,x^\ast;v_1)=\tilde\Lambda(x,x^\ast;v_2)$. Checking carefully the proof of \cite[Theorem~3]{GfrOut14a} allows us to observe that the aforementioned relaxed metric regularity assumption on $M_{\tilde q}$ therein can be replaced by the validity of formula \eqref{EqTanConeGrNormalConeInequ} verified above. Hence we can apply \cite[Theorem~3]{GfrOut14a} to derive the tangent cone and regular normal cone representations
\begin{eqnarray*}
T_{{\rm gph}\,\Hat\partial\delta_\Gamma}(x,x^\ast)=\big\{(v,v^\ast)\big|\;v^\ast\in\nabla^2\la\tilde\lambda,\tilde q)\ra(x)v+N_{K(x,x^\ast)}(v)\big\},
\end{eqnarray*}
\begin{eqnarray*}
\widehat N_{{\rm gph}\,\Hat\partial\delta_\Gamma}(x,x^\ast)=\big\{(w^\ast,w)\big|\;w\in K(x,x^\ast),\;w^\ast\in-\nabla^2\la\tilde\lambda,\tilde q\ra(x)w+ K(x,x^\ast)^*\big\},
\end{eqnarray*}
where $\tilde\lambda$ is an arbitrarily fixed multiplier from $\tilde\Lambda(x,x^\ast;v)$ for some $0\not=v\in K(x,x^\ast)$ if $K(x,x^\ast)\not=\{0\}$ and from $\tilde\Lambda(x,x^\ast)=\tilde\Lambda(x,x^\ast;0)$ otherwise. The obtained formulas easily yield the claimed representations \eqref{EqTangConeConstLambda} and \eqref{EqNormalConeConstLambda} and thus complete the proof of the theorem. $\h$\vspace*{0.05in}

It is worth mentioning that the trivial replacement of an equality by two inequalities as in \eqref{EqIneqSystem} usually does not provide valuable results. However, the imposed MSCQ and BEPP qualification conditions are so powerful, while being fairly nonrestrictive, that they allow us to do it as shown above in the proof of Theorem~\ref{ThTnaConeGrNormalCone}.

The next result is a consequence of Theorem~\ref{ThTnaConeGrNormalCone}(i), which gives us an explicit estimate of the {\em combined second-order subdifferential} \eqref{2nd2} of the constraint indicator function, which is very instrumental in deriving efficient conditions for tilt-stable minimizers in NLPs; see Sections~6 and 7.

\begin{Corollary}[\bf combined second-order subdifferential of the constraint indicator function]\label{LemSecOrd} Assume that both MSCQ and BEPP are satisfied at $\xb$. Then there is a neighborhood $U$ of $\xb$ such that for every $(x,x^\ast)\in\Gr\partial\delta_\Gamma$ with $x\in U$ the following assertion holds: Given any pair $(w,w^\ast)$ with $w^\ast\in\breve\partial^2\delta_\Gamma(x,x^\ast)(w)$, we have $-w\in K(x,x^\ast)$ and
\begin{equation}\label{EqBasicIneqTiltStab}
\skalp{w^\ast,w}\ge\big\la w,\nabla^2\la\lambda,q\ra(x)w\big\ra\;\mbox{ whenever }\;\lambda\in\Lambda(x,x^\ast;-w).
\end{equation}
\end{Corollary}
{\bf Proof.} Let $U$ be neighborhood of $\xb$ such that both MSCQ and BEPP hold for every $x\in\Gamma\cap U$. Fix $x\in\Gamma\cap U$, $x^\ast\in\partial \delta_\Gamma(x)$ and $(w,w^\ast)$ with $w^\ast\in\hat D^\ast\partial \delta_\Gamma(x,x^\ast)(w)$. By definition \eqref{EqRegCoderiv} of the regular coderivative
and representation \eqref{fn2} of the regular normal cone we have
\begin{eqnarray*}
(w^\ast,-w)\in\hat N_{{\rm gph}\,\partial\delta_\Gamma}(x,x^\ast)=\big(T_{{\rm gph}\,\partial\delta_\Gamma}(x,x^\ast)\big)^*.
\end{eqnarray*}
It follows from Theorem~\ref{ThTnaConeGrNormalCone} and Proposition~\ref{LemProxRegDelta} that
$$
\{0\}\times K(x,x^\ast)^*\subset T_{{\rm gph}\,\hat\partial\delta_\Gamma}(x,x^\ast)=T_{{\rm gph}\,\partial\delta_\Gamma}(x,x^\ast).
$$
This implies consequently the relationships
\begin{eqnarray*}
\skalp{w^\ast,0}+\skalp{-w,v^\ast}\le 0\;\mbox{ for all }\;v^\ast\in K(x,x^\ast)^*
\end{eqnarray*}
and hence $-w\in K(x,x^\ast)$. Fixing now any vector $\lambda\in\Lambda(x,x^\ast;-w)$ and using \eqref{EqTanConeGrNormalCone} give us the inclusion $(-w,-\nabla^2\la\lambda,q\ra(x)w)\in T_{{\rm gph}\,\partial\delta_\Gamma}(x,x^\ast)$, and so \eqref{EqBasicIneqTiltStab} is implied by
\begin{eqnarray*}
0\ge\skalp{w^\ast,-w}+\skalp{-w,-\nabla^2\la\lambda,q\ra(x)w}=-\skalp{w^\ast,w}+\big\la w,\nabla^2\la\lambda,q\ra(x)w\big\ra,
\end{eqnarray*}
which completes the proof of the corollary. $\h$\vspace*{0.05in}

The next proposition shows that the stronger CRCQ property yields the additional assumption in Theorem~\ref{ThTnaConeGrNormalCone}(ii) and thus justifies the fulfillment of the simplified formulas \eqref{EqTangConeConstLambda} and \eqref{EqNormalConeConstLambda} therein.

\begin{Proposition}[\bf calculating tangent and regular normal cones under CRCQ]\label{LemMultiplCRCQ} Let CRCQ hold at $\ox\in\Gamma$. Then there is a neighborhood $U$ of $\ox$ such that for every $x\in\Gamma\cap U$ and $x^\ast,v\in\R^n$ satisfying
\begin{eqnarray}\label{I++}
\big\la\nabla q_i(x),v\big\ra=0\;\mbox{ whenever }\;i\in E\cup I^+\;\mbox{ with }\;I^+:=\bigcup_{\lambda\in\Lambda(x,x^\ast)}I^+(\lambda)
\end{eqnarray}
the form  $\la v,\nabla^2\la\cdot,q\ra(x)v\ra$ is constant on $\Lambda(x,x^\ast)$. In particular, we have $\Lambda(x,x^\ast;v)=\Lambda(x,x^\ast)$, and therefore representations \eqref{EqTangConeConstLambda} and \eqref{EqNormalConeConstLambda} are satisfied.
\end{Proposition}
{\bf Proof.} The robustness of CRCQ allows us to proceed in what follows for any $x\in\Gamma$ from some neighborhood $U$ of $\ox$. Consider the case of $\Lambda(x,x^\ast)\ne\emp$ (otherwise the assertion is trivial), fix any $v\in\R^n$ with $\la\nabla q_i(x),v\ra=0$, $i\in E\cup I^+$, and choose the maximal subset $J$ of $E\cup I^+$ such that the gradients $\{\nabla q_i(x)|\;i\in J\}$ are linearly independent. Consider the equations
\begin{equation}\label{EqPartSyst}
q_i\big(x+tv+A^\ast z(t)\big)=0,\quad i\in J,
\end{equation}
where the rows of the $\vert J\vert\times n$ matrix $A$ are given by the gradients $\nabla q_i(x)$, $i\in J$, and where the vectors $z(t)\in\R^{\vert J\vert}$ for any fixed $t\in\R$ are unknown. At $t=0$ we have the trivial solution $z(0)=0$ to \eqref{EqPartSyst} while the Jacobian matrix of this system with respect to $z$ at $t=0$ is the $\vert J\vert\times\vert J\vert$ matrix $AA^\ast$, which is {\em invertible} since the rows of $A$ are linearly independent. Applying the classical {\em implicit function theorem} ensures the existence of $\bar t>0$ and a $C^1$-smooth function $z:(-\bar t,\bar t)\to\R^n$ satisfying the conditions
\[
z(0)=0,\;q_i\big(x+tv+A^\ast z(t)\big)=0\;\mbox{ for all }\;i\in J,\;t\in(-\bar t,\bar t).
\]
By setting $\tilde x(t):=x+tv+A^\ast z(t)$ and differentiating the system \eqref{EqPartSyst} with respect to $t$ we obtain
\[
0=\frac{d}{dt}q_i\big(\tilde x(t)\big)\Big\vert_{t=0}=\big\la\nabla q_i(x),v+A^\ast\frac d{dt}z(0)\big\ra=\big\la\nabla q_i(x),A^\ast\frac d{dt}z(0)\big\ra,\;i\in J,
\]
showing that $AA^\ast\frac d{dt}z(0)=0$ and therefore $\frac d{dt}z(0)=0$. Thus we arrive at the conditions
\begin{eqnarray*}
\tilde x(0)=x,\;\frac{d}{dt}\tilde x(0)=v,\;\mbox{ and }\;q_i\big(\tilde x(t)\big)=0\;\mbox{ for all }\;i\in J,\;t\in(-\bar t,\bar t).
\end{eqnarray*}
It follows from CRCQ that when $t\in(-\bar t,\bar t)$ is sufficiently small, the index set $J$ is the maximal subset of $E\cup I^+$ such that the gradients $\{\nabla q_i(\tilde x(t))|\;i\in J\}$ are linearly independent. Hence for every index $i\in(E\cup I^+)\setminus J$ and small $t\in(-\bar t,\bar t)$ the gradient $\nabla q_i(\tilde x(t))$ can be represented as some linear combination $\sum_{j\in J}\eta_{ij}(t)\nabla q_j(\tilde x(t))$ of $\nabla q_j(\tilde x(t))$, $j\in J$. Employing the standard chain rule tells us that
$$
\frac{d}{dt}q_i\big(\tilde x(t)\big)=\Big\la\nabla q_i\big(x(t)\big),\frac d{dt}\tilde x(t)\Big\ra=\sum_{j\in J}\eta_{ij}(t)\Big\la\nabla q_j\big(\tilde x(t)\big),\frac d{dt}\tilde x(t)\Big\ra=0,\;i\in(E\cup I^+)\setminus J,
$$
and consequently that $q_i(\tilde x(t))=0$ for all $i\in E\cup I^+$ and small $t\in(-\bar t,\bar t)$. Since we also have $\lambda^{(1)}_i=\lambda^{(2)}_i=0$ as $i\in I\setminus\in I^+$ for any $\lambda^{(1)},\lambda^{(2)}\in\Lambda(x,x^\ast)$, it follows that $\la\lambda^{(1)}-\lambda^{(2)},q(\tilde x(t))\ra=0$ if $t\in(-\bar t,\bar t)$ is small enough. Thus by taking into account that $\nabla q(x)^*\lambda^{(1)}=\nabla q(x)^*\lambda^{(2)}=x^\ast$ we get
\begin{eqnarray*}
0&=&\lim_{t\to 0}\frac{\big\la\lambda^{(1)}-\lambda^{(2)},q\big(\tilde x(t)\big)-q(x)\big\ra}{t^2}\\
&=&\lim_{t\to 0}\frac{\big\la\lambda^{(1)}-\lambda^{(2)},\nabla q(x)\big(\tilde x(t)\big)-x\big)\big\ra+\frac 12\big\la\tilde x(t)-x,\nabla^2 \la\lambda^{(1)}-\lambda^{(2)},q\ra(x)\big(\tilde x(t)-x\big)\big\ra+\oo\big(\norm{\tilde x(t)-x}^2\big)}{t^2}\\
&=&\frac 12\big\la v,\nabla^2\la\lambda^{(1)}-\lambda^{(2)},q\ra(x)v\big\ra,
\end{eqnarray*}
which shows that the form $\la v,\nabla^2\la\cdot,q\ra(x)v\ra$ is constant on $\Lambda(x,x^\ast)$ and hence $\Lambda(x,x^\ast;v)=\Lambda(x,x^\ast)$. Since every critical direction $v\in K(x,x^\ast)$ fulfills \eqref{I++} by Proposition~\ref{LemElemProp}, the validity of the claimed representations \eqref{EqTangConeConstLambda} and \eqref{EqNormalConeConstLambda} follows. $\h$

\section{Pointbased Second-Order Sufficient Conditions for Tilt Stability}
\sce

Consider an {\em NLP problem} of minimizing a $C^2$-smooth function $\ph\colon\R^n\to\R$ subject to the constraint system \eqref{EqEqualIneqSystem}, where the equality and inequality constraints are described by $C^2$-smooth functions:
\begin{eqnarray}\label{EqNonlProgr}
\left\{\begin{array}{ll}
\mbox{minimize }\;\varphi(x)\;\mbox{ subject to}\\
q_i(x)=0\;\mbox{ for }\;i\in E\;\mbox{ and }\;q_i(x)\le 0\;\mbox{ for }\;i\in I.
\end{array}\right.
\end{eqnarray}
Using the notation of the previous section, rewrite \eqref{EqNonlProgr} in the unconstrained format
\begin{equation}\label{EqUnconstrProbl}
\mbox{minimize }\;f(x)\;\mbox{ on }\;\R^n,\;\mbox{ where }\;f(x):=\varphi(x)+\delta_\Gamma(x)
\end{equation}
is an extended-real-valued objective. Applying Definition~\ref{DefTiltStab} to the unconstrained problem \eqref{EqUnconstrProbl}, we arrive at the notions of a {\em tilt-stable minimizer} $\ox$ for \eqref{EqNonlProgr}, its {\em modulus} $\kk$, and the {\em exact bound of tilt stability} $\tilt(\ph,q,\ox)$ of the nonlinear program \eqref{EqNonlProgr} at its tilt-stable minimizer $\ox$.

It immediately follows from the subdifferential and coderivative sum rules given in \cite[Proposition~1.107 and Theorem~1.62]{Mo06a}, respectively, that
\[
\hat\partial f(x)=\nabla\varphi(x)+\hat\partial\delta_\Gamma(x),\quad\partial f(x)=\nabla\varphi(x)+\partial\delta_\Gamma(x),\;\mbox{ and}
\]
\[
\breve\partial^2f(x,x^\ast)(v)=\nabla^2\varphi(x)v+\big(\Hat D^*\partial\delta_{\Gamma}\big)\big(x,x^\ast-\nabla\varphi(x)\big)(v)\;\mbox{ whenever }\;x\in\Gamma,\;x^*\in\partial f(x),\;v\in\R^n
\]
for the first-order and second-order subdifferential constructions in \eqref{rs}, \eqref{ls}, and \eqref{2nd2}. Furthermore, we deduce from the definitions of prox-regularity and subdifferential continuity due to Proposition~\ref{LemProxRegDelta} that $f$ in \eqref{EqUnconstrProbl} possesses these properties at any $x\in\Gamma$ close to $\ox$ for $x^\ast\in\partial f(x)$ if both MSCQ and BEPP qualification conditions are satisfied at $\ox$.

By the elementary Fermat rule and sum rule for $\Hat\partial f$ given above we obviously have that any local minimizer $\ox$ for \eqref{EqNonlProgr} fulfills the first-order necessary optimality condition
\begin{eqnarray*}
0\in\hat\partial f(\xb)=\nabla\varphi(\xb)+\hat\partial\delta_\Gamma(\xb),
\end{eqnarray*}
which can be equivalently written (provided that the GCQ $\widehat N_\Gamma(\xb)=\Tlin(\xb)^*$ holds, which is surely the case by \eqref{aba} when MSCQ is satisfied at $\ox$) either as $\Lambda(\xb,-\nabla\varphi(\xb))\ne\emp$ for the set of Lagrange multipliers \eqref{Lambda}, or---more explicitly---in terms of the {\em KKT system}
\begin{eqnarray}\label{EqFO}
0=\nabla_x\Lag(\xb,\lambda)\;\mbox{ for some }\;\lambda\in N_\Theta\big(q(\xb)\big)
\end{eqnarray}
via the classical {\em Lagrange function} defined by
\begin{eqnarray}\label{Lagr}
\Lag(x,\lambda):=\varphi(x)+\la\lambda,q\ra(x)\;\mbox{ for }\;x\in\R^n,\;\lm\in\R^l.
\end{eqnarray}

To formulate our results on tilt stability, define the set of {\em extreme multipliers in critical directions}
\begin{eqnarray}\label{ex-Lambda}
\Lambda_{\cal E}(x,x^\ast;v):=\Lambda(x,x^\ast;v)\cap\E(x,x^\ast)\;\mbox{ for all }\;(x,x^\ast)\in\Gr\hat N_\Gamma,\;v\in K(x,x^\ast),
\end{eqnarray}
which is the collection of extreme points of the multiplier set $\Lambda(x,x^\ast)$ solving the {\em linear program} \eqref{lp}; see the above constructions of ${\cal E}(x,x^*)$ and $\Lambda(x,x^*;v))$. It is well known in linear programming that
$\Lambda_{\cal E}(x,x^\ast;v)\ne\emp$ if and only if both sets $\Lambda(x,x^\ast;v)$ and $\E(x,x^\ast)$ are nonempty; in this case  the set $\Lambda_{\cal E}(x,x^\ast;v)$ precisely reduces to all the extreme points of the convex polyhedron $\Lambda(x,x^\ast;v)$ in \eqref{crit-mult}.\vspace*{0.05in}

Now we are ready to establish the major second-order {\em sufficient} condition for tilt stability in \eqref{EqNonlProgr}, with a prescribed modulus $\kk>0$ and a constructive lower estimate for the exact bound of tilt stability $\tilt(\ph,q,\ox)$, formulated {\em at} the reference point $\ox\in\Gamma$. As the reader can see, this pointbased condition is expressed via the {\em strong positive-definiteness} of the Hessian of the Lagrange function \eqref{Lagr} on the subspace orthogonal to the gradients $\nabla q_i(\ox)$ for the equality and strict complementarity constraint indexes \eqref{I+} generated by {\em extreme multipliers} in all the critical directions \eqref{ex-Lambda} at $(\ox,-\nabla\ph(\ox))$.

\begin{Theorem}[\bf pointbased sufficient condition for tilt-stable minimizers in NLPs with prescribed moduli]\label{ThTiltStabSuff} Given a feasible point $\ox\in\Gamma$ and a number $\kk>0$, suppose that MSCQ, BEPP, and the first-order necessary optimality condition \eqref{EqFO} are satisfied at $\ox$ and that the second-order condition
\begin{eqnarray}\label{2nd-suf1}
\big\la w,\nabla_x^2\Lag(\xb,\lambda)w\big\ra>\frac1\kappa\norm{w}^2\;\mbox{ whenever }\;w\ne 0\;\mbox{ with }\;\la\nabla q_i(\xb),w\ra=0,\;i\in E\cup I^+(\lambda)
\end{eqnarray}
holds for all the $($finitely many$)$ extreme Lagrange multipliers in critical directions
\begin{eqnarray}\label{2nd-suf2}
\lambda\in\bar\Lambda_{\cal E}:=\disp\bigcup_{0\not=v\in K(\xb,-\nabla\varphi(\xb))}\Lambda_{\cal E}\big(\xb,-\nabla\varphi(\xb);v\big).
\end{eqnarray}
Then $\xb$ is a tilt-stable local minimizer for \eqref{EqNonlProgr} with modulus $\kappa$. Furthermore, we have the estimate
\begin{eqnarray}\label{2nd-bound1}
\tilt(\ph,q,\ox)\ge\disp\sup\left\{\frac{\norm{w}^2}{\la w,\nabla^2_x\Lag(\xb,\lambda)w\ra}\Big|\;\lambda\in\bar\Lambda_{\cal E},\;
\la\nabla q_i(\xb),w\ra=0,\;i\in E\cup I^+(\lambda)\right\}\ge 0
\end{eqnarray}
of the exact tilt stability bound of \eqref{EqNonlProgr} at $\ox$ with the convention that $0/0:=0$ in \eqref{2nd-bound1}.
\end{Theorem}
{\bf Proof.} Employing Theorem~\ref{ThTiltStab} and Proposition~\ref{LemProxRegDelta}, it suffices to show that the second-order condition \eqref{2nd-suf1} with $\lm$ from \eqref{2nd-suf2} implies the validity of \eqref{tilt1} for the function $f$ defined in \eqref{EqUnconstrProbl}. Then the exact bound lower estimate \eqref{2nd-bound1} follows directly from \eqref{tilt2} and \eqref{2nd-suf1}.

Suppose on the contrary that \eqref{tilt1} fails while \eqref{2nd-suf2} holds and then find sequences $x_k\to\xb$ and $x_k^\ast\to 0$ as $k\to\infty$ as well as $(w_k,w_k^\ast)\in\gph\breve\partial^2f(x_k,x_k^\ast)$ such that
\begin{equation}\label{EqNegTiltStab}
\skalp{w_k^\ast,w_k}<\frac 1\kappa\norm{w_k}^2\;\mbox{ for all large }\;k\in\N.
\end{equation}
Since $(w_k^\ast,-w_k)\in\hat N_{{\rm gph}\,\partial\delta f}(x_k,x_k^\ast)$ and $w_k\not=0$ by \eqref{EqNegTiltStab}, we may assume that $\norm{w_k}=1$ for all $k$ and select a subsequence $w_k\to w$ with some $w$ from the unit sphere of $\R^n$. It follows from Corollary~\ref{LemSecOrd} that $-w_k\in K(x_k,y_k^\ast)$ with $y_k^\ast:=x_k^\ast-\nabla\varphi(x_k)$. Further, we have by Proposition~\ref{LemElemProp}(iii) that $\Lambda(x_k,y_k^\ast;-w_k)\ne\emp$, and thus the set $\Lambda_{\cal E}(x_k,y_k^\ast;-w_k)$ is also nonempty for each $k\in\N$. Since $w_k^\ast-\nabla^2\varphi(x_k)w_k\in(\hat D^\ast\partial\delta_\Gamma)(x_k,y_k^\ast)(w_k)$ by the above constructions and definition \eqref{2nd2} of the combined second-order subdifferential, we get from the crucial conclusion \eqref{EqBasicIneqTiltStab} of Corollary~\ref{LemSecOrd} that there is a sequence of $\lambda^k\in\Lambda_{\cal E}(x_k,y_k^\ast;-w_k)$ satisfying the inequality
\begin{eqnarray*}
\skalp{w_k^\ast-\nabla^2\varphi(x_k)w_k,w_k}\ge\big\la w_k,\nabla^2\la\lambda^k,q\ra(x_k)w_k\big\ra,
\end{eqnarray*}
which can be rewritten in terms of the Lagrange function \eqref{Lagr} as
\begin{equation}\label{EqNegHessian_w_k}
\skalp{w_k^\ast,w_k}\ge\big\la w_k,\nabla_x^2\Lag(x_k,\lambda^k)w_k\big\ra\;\mbox{ with some }\;\lambda^k\in\Lambda_{\cal E}(x_k,y_k^\ast;-w_k),\quad k\in\N.
\end{equation}
The imposed BEPP at $\ox$ ensures that the sequence $\{\lm^k\}$ from \eqref{EqNegHessian_w_k} is bounded, and hence we find $\bar\lm$ so that $\lm^k\to\bar\lm$ for all $k\to\infty$ without loss of generality. It is easy to see that $\bar\lambda\in N_\Theta(q(\xb))$ and that
\[
\nabla q(\xb)^*\bar\lambda=\lim_{k\to\infty}\nabla q(x_k)^*\lambda^{k}=\lim_{k\to\infty}y_k^\ast=-\nabla\varphi(\xb)
\]
telling us that $\bar\lambda\in\Lambda(\xb,-\nabla\varphi(\xb))$. Let us show next that $\bar\lambda\in\co\E(\xb,-\nabla\varphi(\xb))$.

Assuming the contrary gives us $\bar\lambda=\lambda^e+\lambda^r$ with $\lambda^e\in\co\E(\xb,-\nabla\varphi(\xb))$ and $\lambda^r\not=0$ belonging to the {\em recession cone} of $\Lambda(\xb,-\nabla\varphi(\xb))$, i.e., $\lambda^r_i\ge 0$ for all $i\in I$ and $\nabla q(\xb)^*\lambda^r=0$.
Since
$$
I^+(\lambda^r)\subset I^+(\bar\lambda)\subset I^+(\lambda^k)\subset\I(x_k)
$$
for the index sets \eqref{active} and \eqref{I+} when $k$ is large and since the gradient family $(\nabla q_i(x_k))_{i\in E\cup I^+(\lambda^k)}$ is linearly independent, we have $\lambda^r\in\E(x_k,\nabla q(x_k)^*\lambda^r)$ for such $k$. This clearly contradicts BEPP by
\[
\lim_{k\to\infty}\frac{\norm{\nabla q(x_k)^*\lambda^r}}{\norm{\lambda^r}}=0,
\]
and hence the claimed inclusion $\bar\lambda\in\co\E(\xb,-\nabla\varphi(\xb))$ is verified.

Furthermore, due to the inclusions $I^+(\bar\lambda)\subset I^+(\lambda^k)$ and $-w_k\in K(x_k,y_k^\ast)$ for large $k$, we get $\la\nabla q_i(x_k),-w_k\ra=0$ whenever $i\in E\cup I^+(\lambda^k)$ and thus conclude that $\la\nabla q_i(\xb),w\ra=0$ for all indexes $i\in E\cup I^+(\bar\lambda)$ by passing to the limit as $k\to\infty$. Combining \eqref{EqNegTiltStab} and \eqref{EqNegHessian_w_k} gives us
\begin{equation}\label{EqNegHessian}
\big\la w,\nabla^2\Lag(\xb,\bar\lambda)w\big\ra\le\frac 1\kappa\norm{w}^2
\end{equation}
by the limiting procedure with the limit pair $(w,\bar\lm)$ constructed above.

Consider now the following two cases, which completely cover the situation. In the {\em first case} suppose that $x_k\not=\xb$ for infinitely many $k$ and get $(x_k-\xb)/\norm{x_k-\xb}\to v$ by passing to a subsequence if necessary. Taking into account that $I^+(\bar\lambda)\subset I^+(\lambda^k)\subset\I(x_k)$ for large $k$, we get the relationships
\[
\big\la\nabla q_i(\xb),v\big\ra=\lim_{k\to\infty}\frac{q_i(x_k)-q_i(\xb)}{\norm{x_k-\xb}}\begin{cases}=0&\mbox{if }\;i\in E\cup I^+(\bar\lambda),\\
\le 0&\mbox{if }\;i\in\I(\xb)\setminus I^+(\bar\lambda),\end{cases}
\]
which show that $\nabla q(\xb)^*v\in T_\Theta(q(\xb))$, $0=\la\bar\lambda,\nabla q(\xb)v\ra=-\la\nabla\varphi(\xb),v\ra$, and so  $v\in K(\xb,-\nabla\varphi(\xb))$ by \eqref{crit-cone}. Moreover, $\la\bar\lambda,q(x_k)\ra=0$ when $k$ is sufficiently large, and the conditions $\la\lambda-\bar\lambda,\nabla q(\xb)\ra=0$, $\la\lambda,q(x_k)\ra\le 0$ hold whenever $\lambda\in\Lambda(\xb,-\nabla\varphi(\xb))$. Hence we have
\begin{eqnarray*}
0&\ge&2\limsup_{k\to\infty}\frac{\big\la\lambda-\bar\lambda,q(x_k)\big\ra}{\norm{x_k-\xb}^2}=2\limsup_{k\to\infty}\frac{\big\la\lambda-\bar\lambda,q(x_k)\ra-
\big\la\lambda-\bar\lambda,q(\xb)\big\ra}{\norm{x_k-\xb}^2}\\
&=&\limsup_{k\to\infty}\frac{\big\la x_k-\xb,\nabla^2\la\lambda-\bar\lambda,q\ra(\xb)(x_k-\xb)\big\ra}{\norm{x_k-\xb}^2}=\big\la v,\nabla^2\la\lambda-\bar\lambda,q\ra(\xb)v\big\ra
\end{eqnarray*}
showing that $\bar\lambda\in\Lambda(\xb,-\varphi(\xb);v)$. Since $\bar\lambda\in\co\E(\xb,-\nabla\varphi(\xb))$ by the above, it has the representation
\begin{eqnarray}\label{mu}
\bar\lambda=\disp\sum_{j=1}^N\beta^j\mu^j\;\mbox{ with }\;\mu^j\in\E\big(\xb,-\nabla\varphi(\xb)\big),\;\bb^j>0,\;\sum_{j=1}^N\bb^j=1
\end{eqnarray}
for some $N\in\N$. Taking into account that $\big\la v,\nabla^2\la\mu^j,q\ra(\xb)v\big\ra\le\big\la v,\nabla^2\la\bar\lambda,q\ra(\xb)v\big\ra$ due to the definition of $\Lambda(\xb,-\varphi(\xb);v)$ in \eqref{crit-mult} and that
\[
0=\sum_{j=1}^N\beta^j\big\la v,\nabla^2\la\mu^j,q\ra(\xb)v\big\ra-\big\la v,\nabla^2\la\bar\lambda,q\ra(\xb)v\big\ra,
\]
we conclude that the following relationships are satisfied:
$$
\big\la v,\nabla^2\la\mu^j,q\ra(\xb)v\big\ra=\big\la v,\nabla^2\la\bar\lambda,q\ra(\xb)v\big\ra,\;\mbox{ and so }\;\mu^j\in\Lambda_{\cal E}\big(\xb,-\varphi(\xb);v\big)\;\mbox{ for }\;j=1,\ldots,N.
$$
The latter allows us to use the assumed second-order condition \eqref{2nd-suf1} for $\mu^j$, which implies that
$\la w,\nabla_x^2\Lag(\xb,\bar\lambda)w\ra>\frac 1\kappa\norm{w}^2$ by $I^+(\mu^j)\subset I^+(\bar\lm)$ and hence contradicts \eqref{EqNegHessian}. This justifies the statement of the theorem in the first case under consideration.

In the {\em second case} we have $x_k\not=\xb$ only for finitely many $k$ and so can suppose that $x_k=\xb$ for all $k\in\N$. Since $-w^k\in K(\xb,y_k^\ast)$ as shown above, we easily get that $-w\in K(\xb,-\nabla\varphi(\xb))$ for the limit point $w$. It follows now from \cite[Theorem~5.3.2(2)]{BaGuKlKuTa82} that $\bar\lambda\in\Lambda(\xb,-\nabla\varphi(\xb);-w)$. Representing $\bar\lambda\in\co\E(\xb,-\nabla\varphi(\xb))$  as in \eqref{mu} and using the same arguments as in the first case above, we arrive at a contradiction with \eqref{EqNegHessian} and thus complete the proof of theorem.$\h$\vspace*{0.05in}

Theorem~\ref{ThTiltStabSuff} provides a pointbased second-order sufficient condition for tilt-stable local minimizers of NLPs with a {\em prescribed modulus} $\kk>0$ via the {\em strong} (involving the given modulus $\kk$) positive-definiteness of the Hessian $\nabla^2\Lag(\xb,\bar\lambda)$ in \eqref{2nd-suf1} over the subspace therein with $\lm$ from \eqref{2nd-suf2}. A natural question arises about the sufficiency of the {\em positive-definiteness} counterpart of \eqref{2nd-suf1} and \eqref{2nd-suf2} for tilt stability of $\ox$ with no modulus specified, i.e., in the sense of Definition~\ref{DefTiltStab}(i). The validity of this statement can be justified by using the device similar to the proof of Theorem~\ref{ThTiltStabSuff} while applying instead of Theorem~\ref{ThTiltStab} above (taken from \cite[Theorem~3.5]{MoNg14}) the characterization of tilt stability in the sense of Definition~\ref{DefTiltStab}(i) in the unconstrained format of optimization obtained in \cite[Theorem~1.3]{PolRo98} via the positive-definiteness of the basic second-order subdifferential \eqref{2nd1}. However, the desired result can be also deduced directly from Theorem~\ref{ThTiltStabSuff} as in the following corollary.

\begin{Corollary}[\bf pointbased sufficient condition for tilt-stable minimizers in NLPs with no modulus specified]\label{ThTiltStabSuff1} Assume hat MSCQ, BEPP, and the first-order necessary optimality condition \eqref{EqFO} are satisfied at $\ox$ and that the positive-definiteness condition
\begin{eqnarray}\label{2nd-suf3}
\big\la w,\nabla^2_x\Lag(\xb,\lambda)w\big\ra>0\;\mbox{ for all }\;\lambda\in\bar\Lambda_{\cal E},\;\big\la\nabla q_i(\xb),w\big\ra=0,\;i\in E\cup I^+(\lambda),\;w\ne 0
\end{eqnarray}
holds. Then there is $\kk>0$ such that $\xb$ is a tilt-stable local minimizer with modulus $\kk$ for \eqref{EqNonlProgr}.
\end{Corollary}
{\bf Proof.} Since the set of extreme multipliers $\bar\Lambda_{\cal E}$ is {\em finite} as a subset of extreme points of a convex polyhedron, it is possible to conclude that the positive-definiteness condition \eqref{2nd-suf3} implies its strong counterpart \eqref{2nd-suf2} with the same vectors $\lm$. Indeed, the suitable modulus $\kk>0$ can be constructed so that $\kk^{-1}$ is the minimum of the minimal eigenvalues of the matrices $A_i^*\nabla^2\Lag(\xb,\lambda^i)A_i$, where the columns of $A_i$ form an orthonormal basis of the subspace
$$
\big\{w\in\R^n\big|\;\la\nabla q_i(\ox),w\ra=0,\;i\in E\cup I^+(\lm^i)\big\}\;\mbox{ with }\;\lm^i\in\bar\Lambda_{\cal E}.
$$
The reader may proceed with more details if necessary. $\h$\vspace*{0.05in}

Note that the second-order sufficient conditions in both Theorem~\ref{ThTiltStabSuff} and Corollary~~\ref{ThTiltStabSuff1} trivially hold and ensure tilt stability of $\ox$ if $\bar\Lambda_{\cal E}=\emp$, i.e., when $K(\xb,-\nabla\varphi(\xb))=\{0\}$. However, in this case we can make a more precise statement, which corresponds to $\tilt(\ph,q,\ox)=0$ in \eqref{2nd-bound1}.

\begin{Proposition}[\bf tilt stability with zero exact bound]\label{PropKritConeEmpty} Let MSCQ, BEPP, and the first-order necessary optimality condition \eqref{EqFO} hold at $\xb\in\Gamma$. Suppose further that $K(\xb,-\nabla\varphi(\xb))=\{0\}$. Then for all $\gamma>0$ sufficiently small there is a neighborhood $V$ of the origin in $\R^n$ such that
\begin{eqnarray*}
M_\gamma(v^\ast)=\big\{\xb\big\}\;\mbox{ for all }\;v^\ast\in V,
\end{eqnarray*}
where the argminimum mapping $M_\gg$ is defined in \eqref{EqM_gamma} with $f$ given in \eqref{EqUnconstrProbl}.
\end{Proposition}
{\bf Proof.} The negation of this statement gives us a sequence of $\gamma_k\dn 0$ such that for every neighborhood $V$ of $0\in\R^n$ there exists a vector $v^\ast\in V$ with $M_{\gamma_k}(v^\ast)\not=\{\xb\}$. Using Theorem~\ref{ThTiltStabSuff} and passing to a subsequence if necessary, we can assume
that for every fixed $k\in\N$ there is a neighborhood $V_k$ of $0$ on which $M_{\gamma_k}$ is single-valued and Lipschitz continuous with modulus $k^{-1}$.
Hence for each $k$ we find $v_k^\ast\in V_k\cap\frac{\gamma_k}2\B$ and $x_k\not=\xb$ with $M_{\gamma_k}(v_k^\ast)=\{x_k\}$. Then
$$
\norm{x_k-\xb}\le\frac 1k\norm{v_k^\ast-0}\le\frac 1{2k}\gamma_k,\quad k\in\N,
$$
which shows by definition \eqref{rs} of the regular subdifferential that $0\in\hat\partial(f-\la v_k^\ast,\cdot\ra)(x_k)$. By passing to a subsequence again if needed, we get that $(x_k-\xb)/\norm{x_k-\xb}\to v$ for some unit vector $v\in\R^n$ and that MSCQ and BEPP hold at $x_k$. This justifies the KKT form \eqref{EqFO} of the stationary condition $0\in\hat\partial(f-\la v_k^\ast,\cdot\ra)(x_k)$ and also the existence of a convergent sequence $\lambda_k\to\bar\lambda$ with $\lambda_k\in\E(x_k,-\nabla\varphi(x_k)+v_k^\ast)$. Using the same arguments as in the proof of Theorem~\ref{ThTiltStabSuff} yields $\bar\lambda\in\Lambda(\xb,-\nabla\varphi(\xb))$ and
\[
\big\la\nabla q_i(\xb),v\big\ra=\lim_{k\to\infty}\frac{q_i(x_k)-q_i(\xb)}{\norm{x_k-\xb}}\begin{cases}=0&\mbox{ if }\;i\in E\cup I^+(\bar\lambda),\\
\le 0&\mbox{ if }\;i\in\I(\xb)\setminus I^+(\bar\lambda),\end{cases}
\]
which imply that $v\in\Tlin(\xb)$ and $-\la\nabla\varphi(\xb),v\ra=\la\lambda,\nabla q(\xb)v\ra=0$. This brings us to the contradiction $0\not=v\in K(\xb,-\nabla\varphi(\xb)$ and thus completes the proof of the proposition. $\h$\vspace*{0.05in}

\section{Necessary Conditions and Characterizations of Tilt Stability in NLPs}
\sce

We start with establishing the necessity of the major second-order sufficient condition of Theorem~\ref{ThTiltStabSuff} under additional assumptions involving either {\em nondegeneracy}, or the notion of 2-{\em regularity}. The latter notion was initiated (and named) by Tret'yakov \cite{Tr84} in the case of zero Jacobian and then was strongly developed by Avakov \cite{Av85} whom we mainly follow in the next definition. The symbol $\big[\nabla^2g(\xb)v,w\big]$ stands therein for the $s$-vector column with the quadratic form entries $\la\nabla^2g_i(\ox)v,w\ra$, $i=1,\ldots,s$, generated by the Hessians of all the component $g_i$ of the mapping $g\colon\R^m\to\R^s$.

\begin{Definition}[\bf 2-regularity]\label{2reg} Let $g:\R^m\to\R^s$ be twice Fr\'echet differentiable at $\xb\in\R^m$. We say that $g$ is {\sc 2-regular} at the point $\xb$ in the direction $v\in\R^m$ if for any $p\in\R^s$ the system
\begin{equation}\label{Eq2Reg1}
\nabla g(\xb)u+\big[\nabla^2g(\xb)v,w\big]=p,\quad\nabla g(\xb)w=0
\end{equation}
admits a solution $(u,w)\in\R^m\times\R^m$.
\end{Definition}

Note that Avakov \cite{Av85} used this notion only for directions $v$ satisfying the conditions $\nabla g(\xb)v=0$ and $[\nabla^2 g(\xb)v,v]\in\rge\nabla g(\xb)$ for the {\em range} of the derivative/Jacobian operator $\nabla g(\ox)$.\vspace*{0.05in}

Given $\xb\in\Gamma$, fix a tangent direction $v\in\Tlin(\xb)$ from the linearized constraint cone \eqref{lin-tan} and define the subset of the {\em active inequality} constraint indexes \eqref{active} {\em in the direction} $v$ by
\begin{eqnarray}\label{active1}
\I(\xb;v):=\big\{i\in\I(\xb)\big|\;\la\nabla q_i(\xb),v\ra=0\big\}.
\end{eqnarray}
Introduce further the collection of {\em 2-regularity vectors in the direction} $v\in\Tlin(\xb)$ by
\begin{eqnarray}\label{Xi}
\Xi(\xb;v):=\left\{z\in\R^n\Big|\;\big\la\nabla q_i(\xb),z\big\ra+\big\la v,\nabla^2q_i(\xb)v\big\ra
\begin{cases}=0&\;\mbox{ for }\;i\in E,\\
\le 0&\;\mbox{ for }\;i\in\I(\xb)\end{cases}\right\}
\end{eqnarray}
and consider the corresponding collection of {\em active inequality constraint} indexes \eqref{active1} in this {\em direction}
\begin{eqnarray*}
\J(\xb;v):=\Big\{\J\subset\I(\xb;v)\Big|\,\exists\,z\in\Xi(\xb;v)\;\mbox{with}\;\J=\big\{i\in\I(\xb;v)\big|\;\big\la\nabla q_i(\xb),z\big\ra+ \big\la v,\nabla^2 q_i(\xb)v\big\ra=0\Big\}.
\end{eqnarray*}

The next result shows that 2-regularity of the constraint mapping at the reference point in the given tangent direction built upon equality and ``maximal" active inequality constraints implies a certain {\em parametric LICQ} along a feasible curve with the same active constraint indexes. In what follows we understand a {\em maximal element} of a subset $S$ in a partially ordered set in the usual sense of order theory, i.e., as an element of $S$ that is not smaller than any other element in $S$. It is clear that for $S=\J(\xb;v)$ below a maximal element (by inclusion ``$\subset$") always exists if $\J(\xb;v)\ne\emp$, but it may not be unique.

\begin{Lemma}[\bf parametric LICQ from 2-regularity]\label{Prop2LICQ} Fix $\ox\in\Gamma$, a tangent direction $v\in\Tlin(\xb)$ and a maximal element $\Hat{\cal C}$ of the index subset collection ${\cal C}(\ox;v)$ defined above. Then the 2-regularity of the constraint mapping
$(q_i)_{i\in E\cup\hat\J}$ at $\ox$ in the direction $v$ implies that for every subset $\J\subset\hat\J$ there exists a number $\bar\tau>0$ and a mapping $\hat x:[0,\bar\tau]\to\Gamma$ such that
\begin{eqnarray*}
\hat x(0)=\xb,\;\;\I\big(\hat x(\tau)\big)=\J,\;\;\lim_{\tau\dn 0}\frac{\hat x(\tau)-\xb}\tau=v,
\end{eqnarray*}
and LICQ is satisfied at $\hat x(\tau)$ for every $\tau\in(0,\bar\tau)$.
\end{Lemma}
{\bf Proof.} It is done in \cite[Proposition~4]{GfrOut14b} for the case of inequality constraints, but the given proof goes through by replacing each equality by two inequalities as in the proof of Theorem~\ref{ThTnaConeGrNormalCone}. $\h$\vspace*{0.05in}

Now we are ready to establish the {\em no-gap necessity} of the second-order sufficient condition in Theorem~\ref{ThTiltStabSuff} for tilt-stable minimizers with modulus $\kk>0$, where the strict inequality sign ``$>$" in \eqref{2nd-suf1} is replaced by ``$\ge$" under the extra alternative assumptions: either {\em nondegeneracy in critical directions}, or {\em 2-regularity} of the underlying {\em narrow part} of active constraints in critical directions. Since the latter notion has been formulated in Definition~\ref{2reg}, it remains to introduce the former one.

\begin{Definition} [\bf nondegeneracy in critical direction]\label{nongen-cd} We say that a feasible solution $\ox$ to problem \eqref{EqNonlProgr} {\sc nondegenerates in the critical direction} $v$ if the set of multipliers $\Lambda(\xb,-\nabla\varphi(\xb);v)$ from \eqref{crit-mult} at $(\ox,-\nabla\ph(\ox))$ in the direction $v$ is a singleton.
\end{Definition}

It is clear that this notion is a significant relaxation of the standard notion of nondegeneracity in NLP, which means that the whole set of Lagrange multipliers $\Lambda(\xb,-\nabla\varphi(\xb))$ from \eqref{Lambda} is a singleton.

\begin{Theorem}[\bf no-gap necessary condition for tilt stability with prescribed moduli under either nondegeneracy or 2-regularity]\label{PropNecTiltStab2Reg} Let $\xb$ be a tilt-stable local minimizer with modulus $\kappa>0$ for program \eqref{EqNonlProgr}, and let both MSCQ and BEPP hold at $\xb$. Suppose further that for any nonzero critical direction $0\not=v\in K(\xb,-\nabla\varphi(\xb))$ from \eqref{crit-cone} one of the following assumptions is satisfied:

{\bf(a)} either $\ox$ nondegenerates in the critical directions $v$,

{\bf(b)} or for every extreme multiplier $\lambda\in\Lambda_{\cal E}(\xb,-\nabla\varphi(\xb);v)$ from \eqref{ex-Lambda} there exists a maximal element $\hat\J\in\J(\xb;v)$ such that $I^+(\lambda)\subset\hat\J$ for the strict complementarity index set \eqref{I+} and that the narrow active constraint mapping $(q_i)_{i\in E\cup\hat\J}$ is {\sc 2-regular} at $\ox$ in the direction $v$.\\[1ex]
Then we have the pointbased second-order necessary condition for tilt stability
\begin{eqnarray}\label{2nd-nec1}
\big\la w,\nabla^2_x\Lag(\xb,\lambda)w\big\ra\ge\frac 1\kappa\norm{w^2}\;\mbox{ whenever }\;\lambda\in\bar\Lambda_{\cal E},\;\big\la\nabla q_i(\xb),w\big\ra=0,\;i\in E\cup I^+(\lambda)
\end{eqnarray}
with modulus $\kk$ and with the upper estimate of the exact bound of tilt stability of \eqref{EqNonlProgr} at $\ox$ given by
\begin{eqnarray}\label{2nd-bound2}
\tilt(\ph,q,\ox)\le\disp\sup\left\{\frac{\norm{w}^2}{\la w,\nabla^2_x\Lag(\xb,\lambda)w\ra}\Big|\;\lambda\in\bar\Lambda_{\cal E},\;
\la\nabla q_i(\xb),w\ra=0,\;i\in E\cup I^+(\lambda)\right\}<\infty
\end{eqnarray}
under the convention that $0/0:=0$ in \eqref{2nd-bound2}.
\end{Theorem}
{\bf Proof.} Suppose on the contrary that $\xb$ is a tilt-stable local minimizer with modulus $\kappa$ while
\begin{eqnarray*}
\big\la w,\nabla^2_x\Lag(\xb,\lambda)w\big\ra<\frac 1\kappa\norm{w^2}\;\mbox{ for some }\;\lambda\in\bar\Lambda_{\cal E}\;\mbox{ and }\;w\in\R^n\;\mbox{ with }\;\big\la\nabla q_i(\xb),w\big\ra=0,\;i\in E\cup I^+(\lambda),
\end{eqnarray*}
which obviously yields $w\ne 0$. We now show that there exist a number $\bar\tau>0$ and a mapping $\hat x:[0,\bar\tau]\to\Gamma$ such that $\hat x(0)=\xb$, $\I(\hat x(\tau))=I^+(\lambda)$, LICQ is fulfilled at $\hat x(\tau)$ for every $\tau\in (0,\bar\tau)$, and
\begin{eqnarray*}
\lim_{\tau\downarrow 0}\frac {\hat x(\tau)-\xb}\tau=v,
\end{eqnarray*}
where the nonzero critical direction $0\not=v\in K(\xb,-\nabla\varphi(\xb))$ is chosen such that $\lambda\in \Lambda_{\cal E}(\xb,-\nabla\varphi(\xb);v)$. Observe that
under the assumption made in (b) this follows from Lemma~\ref{Prop2LICQ} with ${\cal C}=I^+(\lambda)$. Hence it remains to consider only the case when
\begin{eqnarray*}
\Lambda\big(\xb,-\nabla\varphi(\xb);v\big)=\Lambda_{\cal E}\big(\xb,-\nabla\varphi(\xb);v\big)=\big\{\lambda\big\}.
\end{eqnarray*}
Recall that by the definition $\lambda$, it is a solution to the linear optimization problem \eqref{lp} with $x=\ox$ and $x^*=-\nabla\ph(\ox)$. Then duality theory in linear optimization ensures the existence of a {\em strictly complementary dual solution} to \eqref{lp}, i.e., some $z\in\R^n$ satisfying
\begin{eqnarray}\label{dual0}
\big\la\nabla q_i(\xb),z\big\ra+\big\la v,\nabla^2q_i(\xb)v\big\ra\begin{cases}=0&\;\mbox{ for }\;i\in E\cup I^+(\lambda),\\
<0&\;\mbox{ for }\;i\in\I(\xb)\setminus I^+(\lambda).\end{cases}
\end{eqnarray}
Taking into account that $\la\nabla q_i(\xb),v\ra=0$ for $i\in E\cup I^+(\lambda)$ gives us
\begin{eqnarray*}
q_i\Big(\xb+\tau v+\frac{\tau^2}2z\Big)=q_i(\xb)+\tau\big\la\nabla q_i(\xb),v\big\ra+\frac{\tau^2}2\Big(\big\la\nabla q_i(\xb),z\big\ra+\big\la v,\nabla^2q_i(\xb)v\big\ra\Big)+\oo(\tau^2)=\oo(\tau^2)
\end{eqnarray*}
whenever $i\in E\cup I^+(\lambda)$. Since $\lambda$ is an extreme point of $\Lambda(\xb,-\nabla\varphi(\xb))$, the constraint gradients $\{\nabla q_i(\xb)|\;i\in E\cup I^+(\lambda)\}$ are linearly independent. Applying the Lyusternik-Graves theorem on metric regularity for smooth mappings, we find positive constants $\gamma$ and $\tau$ such that for every $\tau\in[-\bar\tau,\bar\tau]$ there is is curve $\hat x(\tau)$ satisfying the conditions $q_i(\hat x(\tau))=0$ for $i\in E\cup I^+(\lambda)$ and
\begin{eqnarray*}
\Big\|\hat x(\tau)-\Big(\xb+\tau v+\frac{\tau^2}2z\Big)\Big\|\le\gamma\Big\|\Big(q_i\big(\xb+\tau v+\frac{\tau^2}{2}z\big)\Big)_{i\in E\cup I^+(\lambda)}\Big\|=\oo(\tau^2).
\end{eqnarray*}
Suppose without loss of generality that the gradients $\{\nabla q_i(\hat x(\tau))|\;i\in E\cup I^+(\lambda)\}$ are linearly independent for every $\tau\in[-\bar\tau,\bar\tau]$. Then it follows from $\la\nabla q_i(\xb),v\ra\le 0$ with $i\in\I(\xb)\setminus I^+(\lambda)$ that
\begin{eqnarray*}
q_i\big(\hat x(\tau)\big)=q_i\Big(\xb+\tau v+\frac{\tau^2}2z\Big)+\oo(\tau^2)\le\frac{\tau^2}2\Big(\big\la\nabla q_i(\xb),z\big\ra+\big\la v,\nabla^2q_i(\xb)v\big\ra\Big)+\oo(\tau^2)<0
\end{eqnarray*}
for all $i\in\I(\xb)\setminus I^+(\lambda)$ and $\tau>0$ sufficiently small. Since we also have $q_i(\hat x(\tau))<0$ when $i\in I\setminus\I(\xb)$ and $\tau$ is small enough, it gives us the property
$$
\I\big(\hat x(\tau)\big)=I^+(\lm)\;\mbox{ whenever }\;\tau\in(0,\bar\tau),
$$
this verifies the existence of the curve $\hat x(\cdot)$ with the claimed properties.

Now we pick an arbitrary sequence $\tau_k\downarrow 0$ as $k\to\infty$ with $\tau_k<\bar\tau$ for all $k$ and consider the vectors $x_k:=\hat x(\tau_k)$ and
$x_k^\ast=\nabla q(x_k)^*\lambda$. Denote by $w_k$ the unique optimal solution to the quadratic program:
\begin{eqnarray}\label{quad}
\mbox{minimize }\;\norm{u-w}^2\;\mbox{ subject to }\;\big\la\nabla q_i(x_k),u\big\ra=0\;\mbox{ for all }\;i\in E\cup I^+(\lambda).
\end{eqnarray}
Employing standard arguments in such settings (see, e.g., in the proof of \cite[Theorem~8.2]{KlKum02}) shows that $w_k\to w$ as $k\to\infty$. Moreover, it follows from Proposition~\ref{LemElemProp}(i) on the description of the critical cone \eqref{crit-cone} and from the constraint structure in \eqref{quad} that $-w_k\in K(x_k,x_k^\ast)$ and also that
\begin{eqnarray*}
\widehat N_{{\rm gph}\,\partial\delta_\Gamma}(x_k,x_k^\ast)=\widehat N_{{\rm gph}\,\hat\partial\delta_\Gamma}(x_k,x_k^\ast)=\big\{(u^\ast,u)\big|\;
u\in K(x_k,x_k^\ast),\;u^\ast\in-\nabla^2\big\la\lambda,q\big\ra(x_k)u+K(x_k,x_k^\ast)^*\big\}
\end{eqnarray*}
by further applying Theorem~\ref{ThTnaConeGrNormalCone} and Proposition~\ref{LemProxRegDelta}. Therefore we get
\begin{eqnarray*}
\big(\nabla^2\la\lambda,q\ra(x_k)w_k,-w_k\big)\in\widehat N_{{\rm gph}\,\partial\delta_\Gamma}(x_k,x_k^\ast),\;\mbox{ and so }\;\nabla_x^2\Lag(x_k,\lambda)w_k\in\breve\partial^2f\big(x_k,\nabla\varphi(x_k)+x_k^\ast\big)(w_k).
\end{eqnarray*}
Since $\nabla\varphi(x_k)+x_k^\ast=\nabla_x\Lag(x_k,\lambda)\to 0$ as $k\to\infty$, it follows from Theorem~\ref{ThTiltStab} that
\begin{eqnarray*}
\big\la w_k,\nabla_x^2\Lag(x_k,\lambda)w_k\big\ra\ge\frac 1\kappa\norm{w_k}^2\;\mbox{ for all large }\;k\in\N.
\end{eqnarray*}
By passing to the limit as $k\to\infty$, this clearly contradicts the assumption made at the beginning of the proof of this theorem, and hence we arrive at the
necessary condition \eqref{2nd-nec1} for tilt stability. The exact bound estimate \eqref{2nd-bound2} easily follows from  \eqref{2nd-nec1}, and thus we are done. $\h$\vspace*{0.05in}

The next result is a consequence of Theorem~\ref{PropNecTiltStab2Reg} ensuring the necessity of the pointbased {\em positive-definiteness} condition \eqref{2nd-suf3} from Corollary~\ref{ThTiltStabSuff1} for tilt-stable minimizers of \eqref{EqNonlProgr} {\em with no modulus specified} under the mild assumptions of Theorem~\ref{PropNecTiltStab2Reg}.

\begin{Corollary}[\bf pointbased necessary condition for tilt-stable minimizers in NLPs with no modulus specified]\label{nec2reg}  Let $\xb$ be a tilt-stable local minimizer of \eqref{EqNonlProgr} under the assumptions of Theorem~{\rm\ref{PropNecTiltStab2Reg}}. Then the second-order positive-definiteness condition \eqref{2nd-suf3} is satisfied.
\end{Corollary}
{\bf Proof.} If $\ox$ is a tilt-stable minimizers of \eqref{EqNonlProgr}, then by Definition~\ref{DefTiltStab}(i) applied to the function $f$ from \eqref{EqUnconstrProbl} there is $\kk>0$ such that $\ox$ is tilt stable for \eqref{EqNonlProgr} with modulus $\kk$ as formulated in Definition~\ref{DefTiltStab}(ii). Thus we get condition \eqref{2nd-nec1} by Theorem~\ref{PropNecTiltStab2Reg}, which obviously implies \eqref{2nd-suf3}. $\h$\vspace*{0.05in}

Now we are ready to present {\em complete characterizations} of tilt-stable minimizers for \eqref{EqNonlProgr} with and without prescribed moduli, which are combinations of the results obtained above while definitely deserve to be formulated as a theorem. Moreover, the following theorem contains the {\em precise pointbased formula} for calculating the exact bound of tilt stability.

\begin{Theorem}[\bf second-order characterizations of tilt stability for NLPs under either nondegeneracy or 2-regularity]\label{CorEquivTiltStab2reg} Let $\ox\in\Gamma$ be a feasible solution to \eqref{EqNonlProgr} satisfying MSCQ, BEPP, and the first-order optimality condition \eqref{EqFO}. Suppose further that for every $0\not=v\in K(\xb,-\nabla\varphi(\xb))$ either assumptions in {\rm(a)} or in {\rm(b)} of Theorem~{\rm\ref{PropNecTiltStab2Reg}} are also satisfied. Then the following assertions hold:

{\bf (i)} Given $\kk>0$, the point $\ox$ is a tilt-stable minimizer of \eqref{EqNonlProgr} with any modulus $\kk'>\kk$ if and only if the second-order condition \eqref{2nd-suf3} is fulfilled.

{\bf (ii)} The point $\ox$ is tilt-stable minimizer of \eqref{EqNonlProgr} with some modulus $\kk>0$ if and only if we have the positive-definiteness condition over the extreme multipliers formulated in \eqref{2nd-suf3}.

Furthermore, the exact bound of tilt stability of \eqref{EqNonlProgr} at $\ox$ is finite and calculated by
\begin{eqnarray*}
\tilt(\ph,q,\ox)=\disp\sup\left\{\frac{\norm{w}^2}{\la w,\nabla^2_x\Lag(\xb,\lambda)w\ra}\Big|\;\lambda\in\bar\Lambda_{\cal E},\;
\la\nabla q_i(\xb),w\ra=0,\;i\in E\cup I^+(\lambda)\right\},
\end{eqnarray*}
where we use the convention that $0/0:=0$ as above.
\end{Theorem}
{\bf Proof.} It follows from the combination of the results obtained in Theorem~\ref{ThTiltStabSuff} and Theorem~\ref{PropNecTiltStab2Reg} for assertion (i) and in Corollary~\ref{ThTiltStabSuff1} and Corollary~\ref{nec2reg} for assertion (ii). $\h$\vspace*{0.05in}

Note that the second-order necessary conditions for tilt stability obtained above (and hence the characterizations of Theorem~\ref{CorEquivTiltStab2reg}) involves a certain {\em nondegeneracy} in critical directions---either explicitly assumed in Theorem~\ref{PropNecTiltStab2Reg}(a), or via 2-regularity in  Theorem~\ref{PropNecTiltStab2Reg}(b) that reduces to nondegeneracy by Lemma~\ref{Prop2LICQ}. The next result shows that these nondegeneracy assumptions can be avoided if our basic qualification conditions MSCQ and BEPP are replaced by the stronger CRCQ at the reference point. Observe that the pointbased second-order characterizations of tilt stability obtained in the new setting are somewhat different from those in Theorem~\ref{PropNecTiltStab2Reg} and are expressed via the set of {\em all the Lagrange multipliers} \eqref{Lambda}, while still being pointbased and constructive. It is also worth mentioning that, in the absence of LICQ,
the assumptions of Theorem~\ref{PropNecTiltStab2Reg}(b) and Theorem~\ref{ThTiltStabCRCQ} are strictly {\em complementary} to each other. Indeed, the assumptions of Theorem~\ref{PropNecTiltStab2Reg}(b) imply that the gradients of the active inequality constraints are {\em linearly independent}, while CRCQ imposed in Theorem~\ref{ThTiltStabCRCQ} requires its {\em linear dependence} around the reference point.\vspace*{0.05in}

In the new theorem  presented below we exclude the case of $K(\xb,-\nabla\varphi(\xb))=\{0\}$, which has been already considered in Proposition~\ref{PropKritConeEmpty}.

\begin{Theorem}[\bf second-order characterizations of tilt stability for NLPs under CRCQ]\label{ThTiltStabCRCQ} Let $\ox\in\Gamma$ be a feasible solution to \eqref{EqNonlProgr} satisfying CRCQ and the first-order optimality condition \eqref{EqFO}. Assume further that the gradients of the equality constraints $\{\nabla q_i(\xb)|\;i\in E\}$ are linearly independent and that $K(\xb,-\nabla\varphi(\xb))\not=\{0\}$ for the critical cone \eqref{crit-cone}. Then the following assertions hold:

{\bf (i)} Given $\kk>0$, $\ox$ is a tilt-stable minimizer of \eqref{EqNonlProgr} with any modulus $\kk'>\kk$ if and only if
\begin{eqnarray}\label{2nd-nec3}
\big\la w,\nabla^2_x\Lag(\xb,\lambda)w\big\ra\ge\frac 1\kappa\norm{w}^2\;\mbox{ for all }\;\lambda\in\Lambda\big(\ox,-\nabla\ph(\ox)\big),\;\big\la\nabla q_i(\xb),w\big\ra=0,\;i\in E\cup I^+,
\end{eqnarray}
where $I^+$ is defined in \eqref{I++} with $x=\ox$ and $x^*=-\nabla\ph(\ox)$ while $I^+(\lm)$ is taken from \eqref{I+}. Moreover, the latter is equivalent to the condition:\\[1ex]
For every $w\in\R^n$ with $\big\la\nabla q_i(\xb),w\big\ra=0$ as $i\in E\cup I^+$ there is $\lambda\in\Lambda(\ox,-\nabla\ph(\ox))$ such that
\begin{eqnarray}\label{2nd-nec3c}
\big\la w,\nabla^2_x\Lag(\xb,\lambda)w\big\ra\ge\frac 1\kappa\norm{w}^2.
\end{eqnarray}

{\bf (ii)} The point $\ox$ is a tilt-stable minimizer of \eqref{EqNonlProgr} without modulus specified if and only if
\begin{equation}\label{2nd-nec3e}
\big\la w,\nabla^2_x\Lag(\xb,\lambda)w\big\ra>0\;\mbox{ for all }\;\lambda\in\Lambda\big(\ox,-\nabla\ph(\ox)\big),\;w\not=0,\;\big\la\nabla q_i(\xb),w\big\ra=0,\;i\in E\cup I^+,
\end{equation}
which is equivalent to positive-definiteness condition: for every $0\not=w\in\R^n$ with $\la\nabla q_i(\xb),w\ra=0$ whenever $i\in E\cup I^+$ there is a multiplier $\lambda\in\Lambda\big(\ox,-\nabla\ph(\ox)\big)$ such that $\big\la w,\nabla^2_x\Lag(\xb,\lambda)w\big\ra>0$.

In any of these cases the exact bound of tilt stability of \eqref{EqNonlProgr} at $\ox$ is finite and calculated by
\begin{eqnarray*}
\tilt(\ph,q,\ox)=\disp\sup\left\{\frac{\norm{w}^2}{\la w,\nabla^2_x\Lag(\xb,\lambda)w\ra}\Big|\;\lambda\in\Lambda\big(\ox,-\nabla\ph(\ox)\big),\;
\la\nabla q_i(\xb),w\ra=0,\;i\in E\cup I^+\right\}
\end{eqnarray*}
with the convention that $0/0:=0$ as above.
\end{Theorem}
{\bf Proof.} First we justify the {\em sufficiency} of \eqref{2nd-nec3} for the tilt stability of $\ox$ with any modulus $\kk'>\kk$. Pick any $\tilde\lambda\in\Lambda(\xb,-\varphi(\xb))$ with $I^+(\tilde\lambda)=I^+$ by Proposition~\ref{LemElemProp}(ii) and proceed similarly to the proof of Theorem~\ref{ThTiltStabSuff}. Suppose on the contrary that there are sequences $x_k\to\xb$ and $x_k^\ast\to 0$ and pairs $(w_k,w^*_k)\in\gph\breve\partial^2f(x_k,x_k^\ast)(w_k)$ satisfying
\begin{eqnarray}\label{suff1}
\big\la w_k^\ast,w_k\big\ra<\frac 1{\kappa'}\norm{w_k}^2\;\mbox{ for some }\;\kk'>\kk\;\mbox{ and all }\;k\in\N.
\end{eqnarray}
Let $y_k^\ast$, $\lambda^k$, $\bar\lambda$, and $w$ be chosen as in the proof of Theorem~\ref{ThTiltStabSuff}. Since $I^+(\bar\lambda)\subset I^+(\lambda^k)\subset \I(x_k)$ for all $k$ sufficiently large and since we have the equality
\begin{equation}\label{EqDegeneracy}
\sum_{i\in E\cup I^+(\bar\lambda)}\big(\tilde\lambda_i-\bar\lambda_i\big)\nabla q_i(\xb)+\sum_{i\in I^+\setminus I^+(\bar\lambda)}\tilde\lambda_i\nabla q_i(\xb)=0
\end{equation}
by the definition of $I^+(\lm)$ in \eqref{I+}, it follows from the result of Lu \cite[Proposition~1]{Lu11} involving CRCQ that $I^+\subset\I(x_k)$, i.e., all the constraints for $i\in I^+$ are active at $x_k$. Hence the critical directions $-w_k\in K(x_k,y_k^\ast)$ satisfy the relationships
\begin{eqnarray*}
\big\la\nabla q_i(x_k),-w_k\big\ra=0\;\mbox{ if }\;i\in E\cup I^+(\lambda^k)\;\mbox{ and }\;\big\la\nabla q_i(x_k),-w_k\big\ra\le 0\;\mbox{ if }\;i\in \I(x_k)\setminus I^+(\lambda^k)
\end{eqnarray*}
for all large $k$, which readily ensure their limiting counterparts
\begin{eqnarray*}
\big\la\nabla q_i(\xb),-w\big\ra=0\;\mbox{ if }\;i\in E\cup I^+(\bar\lambda)\;\mbox{ and }\;\big\la\nabla q_i(\xb),-w\big\ra\le 0\;\mbox{ if }\;i\in I^+\setminus I^+(\bar\lambda).
\end{eqnarray*}
By scalar multiplication of \eqref{EqDegeneracy} and $w$ with taking into account that $\tilde\lambda_i>0$ if $i\in I^+\setminus I^+(\bar\lambda)$ we get
\begin{eqnarray*}
\big\la\nabla q_i(\xb),w\big\ra=0\;\mbox{ for }\;i\in I^+\setminus I^+(\bar\lambda),\;\mbox{ and so }\;\big\la\nabla q_i(\xb),w\big\ra=0\;\mbox{ for }\;i\in E\cup I^+.
\end{eqnarray*}
Proceeding then as in the proof of Theorem~\ref{ThTiltStabSuff} gives us a contradiction with \eqref{suff1} and thus verifies the sufficiency part of this theorem with the lower estimate ``$\ge$" in the exact bound formula.

To prove next the {\em necessity} of \eqref{2nd-nec3} for tilt stability in (i) and hence the upper estimate of $\tilt(\ph,q,\ox)$, suppose on the contrary that $\xb$ is a tilt-stable local minimizer with modulus $\kappa'$ for every $\kappa'>\kappa$, but there are vectors $\lambda\in \Lambda(\xb,-\nabla\varphi(\xb))$ and $w\in\R^n\setminus\{0\}$ satisfying
\begin{eqnarray}\label{suff2}
\big\la w,\nabla^2_x\Lag(\xb,\lambda)w\big\ra<\frac 1\kappa\norm{w^2}\;\mbox{ with }\;\big\la\nabla q_i(\xb),w\big\ra=0\;\mbox{ for all }\;i\in E\cup I^+.
\end{eqnarray}
We can clearly choose $\kappa'>\kappa$ so close to $\kappa$ that inequality \eqref{suff2} holds and can suppose by Proposition~\ref{LemMultiplCRCQ} that $\lambda=\tilde\lambda$ without loss of generality. Proposition~\ref{LemElemProp}(ii) allows us to select a critical direction
$v\in K(\xb,-\nabla\varphi(\xb))$ satisfying the conditions
\begin{eqnarray*}
\big\la\nabla q_i(\xb),v\big\ra=0\;\mbox{ for }\;i\in E\cup I^+\;\mbox{ and }\;\big\la\nabla q_i(\xb),v\big\ra<0\;\mbox{ for }\;i\in\I(\xb)\setminus I^+.
\end{eqnarray*}
Following the proof of Proposition~\ref{LemMultiplCRCQ} under CRCQ, we find a $C^1$-smooth mapping $\tilde x:(-\bar t,\bar t)\to\R^n$ with $\tilde x(0)=\xb$, $\frac d{dt}\tilde x(0)=v$, and $q_i(\tilde x(t))=0$ when $i\in E\cup I^+$ and $t\in(-\bar t,\bar t)$. This yields $\frac d{dt}q_i(\tilde x(0))=\nabla q_i(\xb)v<0$ whenever $i\in\I(\xb)\setminus I^+$, and thus $\tilde x(t)\in\Gamma$ and $\I(\tilde x(t))=I^+$ for all small $t>0$.

To complete the proof of this part, we proceed similarly to the case of Theorem~\ref{PropNecTiltStab2Reg} selecting an arbitrary sequence $\tau_k\downarrow 0$, setting $x_k:=\tilde x(\tau_k)$ and $x_k^\ast:=\nabla q(x_k)^*\lambda$, and denoting by $w_k$ the unique optimal solution to the quadratic optimization program:
\begin{eqnarray*}
\mbox{minimize }\;\norm{u-w}^2\;\mbox{ subject to }\;\big\la\nabla q_i(x_k),u\big\ra=0\;\mbox{ for all }\;i\in J,
\end{eqnarray*}
where $J$ is a maximal subset of $E\cup I^+$ such that the gradient vectors $\{\nabla q_i(\xb)|\;i\in J\}$ is linearly independent. Then $w_k\to w$ as in the proof of Theorem~\ref{PropNecTiltStab2Reg}, while the assumed CRCQ tells us that $\la\nabla q_i(x_k),w_k\ra=0$ for all $i\in E\cup I^+$ and $k$ sufficiently large. Following again the arguments of the latter theorem and taking into account that $\Lambda(\ox,-\nabla\ph(\ox);v)=\Lambda(\ox,-\nabla\ph(\ox))$ by Proposition~\ref{LemMultiplCRCQ} valid under CRCQ, we employ the regular normal cone representation \eqref{EqNormalConeConstLambda} with $\lm$ from \eqref{mult1} and thus arrive at a contradiction with \eqref{suff2}, which establishes the ``if and if" statement in (i). The equivalence of \eqref{2nd-nec3} for the validity of \eqref{2nd-nec3c} with {\em some} $\lambda\in\Lambda(\ox,-\nabla\ph(\ox))$ follows immediately from Proposition~\ref{LemMultiplCRCQ}.\vspace*{0.05in}

Verifying finally assertion (ii), it remains to observe by Proposition~\ref{LemElemProp}(ii) that the positive-definiteness condition \eqref{2nd-nec3e} with the selected $\tilde\lm\in\Lambda(\ox,-\nabla\ph(\ox))$ from that proposition yields
\begin{eqnarray*}
\big\la w,\nabla^2_x\Lag(\xb,\tilde\lambda)w\big\ra\ge\mu\norm{w^2}\;\mbox{ whenever }\;\big\la\nabla q_i(\xb),w\big\ra=0\;\mbox{ as }\;i\in E\cup I^+
\end{eqnarray*}
for this fixed $\tilde\lm$, where the number $\mu$ is positive and is defined by
\begin{eqnarray*}
\mu:=\inf\{\big\la w,\nabla^2_x\Lag(\xb,\tilde\lambda)w\big\ra\big|\;\norm{w}=1,\;\big\la\nabla q_i(\xb),w\big\ra=0,\;i\in E\cup I^+\big\}.
\end{eqnarray*}
This completes the proof of the theorem by taking into account the discussions above. $\h$\vspace*{0.05in}

The final result of this section shows that, as far as second-order analysis is concerned, the pointbased sufficient condition \eqref{2nd-suf3} of Corollary~\ref{ThTiltStabSuff1} is also {\em necessary} for tilt stability in the sense of Definition~\ref{DefTiltStab} {\em without} either nondegeneracy or CRCQ requirements of Theorems~\ref{PropNecTiltStab2Reg} and \ref{ThTiltStabCRCQ}. The only assumption needed for this statement is the pointbased
{\em SOSCMS} property \eqref{SOSCMS} by the first author \cite{Gfr11} that has already been discussed in Section~4. The exact meaning of the theorem below is that violating \eqref{2nd-suf3} at $\ox\in\Gamma$ for the given NLP \eqref{EqNonlProgr} yields the violation of tilt stability in a modified NLP with the same cost function and the same values of the constraint functions and their first and second derivatives at $\ox$ as in \eqref{EqNonlProgr}. Thus \eqref{2nd-suf3} is in fact an {\em unimprovable pointbased characterization} of tilt-stable minimizers for NLPs with $C^2$-smooth data under the mild SOSCMS assumption.

\begin{Theorem}[\bf pointbased second-order characterization of tilt stability under SOSCMS]\label{ThNecCond2} Let $\ox\in\Gamma$ satisfy the first-order optimality condition \eqref{EqFO} in NLP \eqref{EqNonlProgr} as well as SOSCMS in the form
\begin{equation}\label{EqSOSCMS}
\big[\lambda\in N_\Theta\big(q(\xb)\big),\;\nabla q(\xb)^\ast\lambda=0,\;\big\la u,\nabla^2\la\lambda,q\ra(\xb)u\big\ra\ge 0\big]\Longrightarrow\lm=0.
\end{equation}
Suppose further that the negation of \eqref{2nd-suf3} holds, i.e.,
\begin{eqnarray*}
\big\la w,\nabla_x^2\Lag(\xb,\lambda)w\big\ra\le 0\;\mbox{ for some }\;\lm\in\bar\Lambda_{\cal E}\;\mbox{ and }\;w\ne 0\;\mbox{ with }\;\la\nabla q_i(\xb),w\ra=0,\;i\in E\cup I^+(\lambda).
\end{eqnarray*}
Then there exist $C^2$-smooth functions $\hat q_i\colon\R^n\to\R$ as $i=1,\ldots,l$ satisfying
\begin{eqnarray*}
\hat q_i(\xb)=q_i(\xb),\;\nabla\hat q_i(\xb)=\nabla q_i(\xb),\;\mbox{ and }\;\nabla^2\hat q_i(\xb)=\nabla^2q_i(\xb)\;\mbox{ for all }\;i=1,\ldots,l
\end{eqnarray*}
and such that $\xb$ is not a tilt-stable local minimizer of the modified nonlinear program
\begin{equation}\label{EqPertProbl}
\mbox{minimize }\;\ph(x)\;\mbox{ subject to }\;\hat q_i(x)=0\;\mbox{ for }\;i\in E\;\mbox{ and }\;\hat q_i(x)\le 0\;\mbox{ for }\;i\in I.
\end{equation}
\end{Theorem}
{\bf Proof.} Take a critical direction $0\not=v\in K(\xb,-\varphi(\xb))$ for which $\lambda\in\Lambda_{\cal E}(\xb,-\nabla\varphi(\xb);v)$ and suppose without loss of generality that $\norm{v}=1$. Recall by the definition of $\Lambda_{\cal E}(\xb,-\nabla\varphi(\xb);v)$ in \eqref{ex-Lambda} that
$\lm$ solves the linear program \eqref{lp} with $x=\ox$ and $x^*=-\nabla\ph(\ox)$. Consider now the problem
\begin{eqnarray*}
\mbox{minimize }\;\big\la\nabla\varphi(\xb),z\big\ra\;\mbox{ subject to }\;\big\la\nabla q_i(\xb),z\big\ra+\big\la v,\nabla^2q_i(\xb)v\big\ra\begin{cases}=0,&i\in E,\\\le 0,&i\in\I(\xb),\end{cases}
\end{eqnarray*}
which is {\em dual} to \eqref{lp} with $x=\ox$ and $x^*=-\nabla\ph(\ox)$. Since $\lm$ solves \eqref{lp}, classical duality in linear programming ensures the existence of $\tilde z\in\R^n$ that solves the dual program and satisfies
\begin{eqnarray*}
\big\la\nabla q_i(\xb),\tilde z\big\ra+\big\la v,\nabla^2q_i(\xb)v\big\ra\begin{cases}=0,&i\in E\cup I^{+}(\lambda),\\\le 0,&i\in\I(\xb)\setminus I^+(\lambda).\end{cases}
\end{eqnarray*}
 Denoting now $z=\tilde z+\alpha v$ for some $\alpha>0$ sufficiently large, we have
\begin{eqnarray*}
\big\la v,z\big\ra>0\;\mbox{ and }\;\big\la\nabla q_i(\xb),z\big\ra+\big\la v,\nabla^2q_i(\xb)v\big\ra\begin{cases}=0,&i\in E\cup I^{+}(\lambda),\\\le 0,&i\in\I(\xb)\setminus I^+(\lambda).\end{cases}
\end{eqnarray*}
Furthermore, the reader can directly check the following inequalities:
\begin{eqnarray*}
1+2\big\la z,v\ra\la x-\xb,v\ra\ge 1-2\big\la z,v\big\ra\norm{x-\xb}>0\;\mbox{ whenever }\;\|x-\ox\|<2r\;\mbox{ with }\;r:=\frac 1{4\big\la z,v\big\ra},
\end{eqnarray*}
\begin{eqnarray*}
1+2\big\la z,v\ra\la x-\xb,v\ra+64\big\la z,v\big\ra^3\big(\norm{x-\xb}-r\big)^3&\ge&1-2\big\la z,v\big\ra\norm{x-\xb}+64\big\la z,v\big\ra^3r^2\big(\norm{x-\xb}-r\big)\\
&=&2\big\la z,v\big\ra\norm{x-\xb}>0\;\mbox{ whenever }\;\|x-\ox\|\ge 2r.
\end{eqnarray*}
These relationships allow us to define the real-valued function $\vt:\R^n\to\R$ by
\begin{eqnarray}\label{vt}
\vt(x):=\frac{-1+\sqrt{1+2\big\la z,v\ra\la x-\xb,v\ra+64\big\la z,v\big\ra^3\max\big\{\norm{x-\xb}-r,0\big\}^3}}{\big\la z,v\big\ra},
\end{eqnarray}
which is clearly twice continuously differentiable in $\R^n$ satisfying the condition
\begin{eqnarray*}
\vt\Big(\xb+tv+\frac 12t^2 z\Big)=\frac{-1+\sqrt{1+2t\big\la z,v\big\ra+t^2\big\la z,v\big\ra^2}}{\big\la z,v\big\ra}=t\;\mbox{ whenever }\;\Big\|tv+\frac 12 t^2 z\Big\|\le r
\end{eqnarray*}
together with $\nabla\vt(\xb)=v$ and $\nabla^2\vt(\xb)=\la z,v\ra vv^{(*)}$, where $vv^{(*)}$ indicated the matrix multiplication of the vector column $v\in\R^n$ by the vector row $v^{(*)}$. Consider next the index set
\begin{eqnarray*}
\hat I:=\big\{i\in\I(\xb;v)\big|\;\big\la\nabla q_i(\xb),z\big\ra+\big\la v,\nabla^2q_i(\xb)v\big\ra=0\big\}\supset I^+(\lambda)
\end{eqnarray*}
and by using \eqref{vt} define the new constraint functions as follows
\begin{eqnarray*}
\hat q_i(x):=\begin{cases}q_i(x)-q_i\big(\xb+\vt(x)v+\frac 12\vt(x)^2z\big)&\mbox{ for }\;i\in E\cup I^+(\lambda),\\
q_i(x)-q_i\big(\xb+\vt(x)v+\frac 12\vt(x)^2z\big)-\norm{x-\xb}^4&\mbox{ for }\;i\in\hat I\setminus I^+(\lambda),\\
q_i(x)&\mbox{ for }\;i\in I\setminus\hat I.\end{cases}
\end{eqnarray*}
This gives us the following relationships with the original constraint functions:
\begin{eqnarray*}
\nabla\hat q_i(\xb)&=&\nabla q_i(\xb)-\big\la\nabla q_i(\xb),v\big\ra v=\nabla q_i(\xb),\\
\nabla^2\hat q_i(\xb)&=&\nabla^2q_i(\xb)-\big(\big\la v,\nabla^2q_i(\xb)v\big\ra+\big\la\nabla q_i(\xb),z\big\ra-
\la z,v\ra\la\nabla q_i(\xb),v\ra\big)vv^{(*)}=\nabla^2q_i(\xb)
\end{eqnarray*}
whenever $i\in E\cup\hat I$. Furthermore, for $0<t$ with $0<\norm{tv+\frac 12 t^2 z}\le r$ it holds
\[
\hat q_i\Big(\xb+tv+\frac 12 t^2z\Big)=\begin{cases}0&\mbox{ for }\;i\in E\cup I^+(\lambda),\\-\norm{tv+t^2z}^4<0&\mbox{ for }\;i\in\hat I\setminus I^+(\lambda).\end{cases}
\]
Since $I\setminus\hat I$ can be partitioned into the sets $I\setminus\I(\xb)$, $\I(\xb)\setminus\I(\xb;v)$, and $\I(\xb;v)\setminus\hat I$ with
\begin{eqnarray*}
\begin{array}{ll}
q_i(\xb)<0\mbox{ for }\;I\setminus \I(\xb),\quad q_i(\xb)=0,\;\big\la\nabla q_i(\xb),v\big\ra<0\;\mbox{ for }\;i\in\I(\xb)\setminus \I(\xb;v),\mbox{ and }\\
q_i(\xb)=0,\;\big\la\nabla q_i(\xb),v\big\ra=0,\;\big\la\nabla q_i(\xb),z\big\ra+\big\la v,\nabla^2q_i(\xb)v\big\ra< 0\;\mbox{ for }\;i\in\I(\xb;v)\setminus\hat I
\end{array}
\end{eqnarray*}
and by the validity of the representation
\begin{eqnarray*}
\begin{array}{ll}
\disp\hat q_i\Big(\xb+tv+\frac 12 t^2b\Big)&=\disp q_i\Big(\xb+tv+\frac 12 t^2z\Big)=q_i(\xb)+t\big\la\nabla\disp q_i(\xb),v\big\ra+\frac{t^2}2\Big(\big\la\nabla q_i(\xb),z\big\ra\\
&+\big\la v,\nabla^2q_i(\xb)v\big\ra\Big)+\oo(t^2),
\end{array}
\end{eqnarray*}
we conclude that $\hat q_i(\xb+tv+\frac 12 t^2z)<0$ for all $i\in I\setminus\hat I$ and all $t>0$ sufficiently small. It follows from
\eqref{EqSOSCMS} that SOSCMS \eqref{SOSCMS} is satisfied for $\Hat q$ with $E_2=I_2=\emp$, and so Theorem~\ref{ThSSOSCMS} shows that both MSCQ and BEPP holds at $\ox$ for the modified constraint system $\hat q(x)\in\Theta$ with $\Th$ from \eqref{Gamma}.

To complete the proof, pick an arbitrary sequence $t_k\downarrow 0$ as $k\to\infty$, denote
$$
x_k:=\xb+t_kv+\disp\frac12t_k^2z\;\mbox{ and }\;x_k^\ast:=\nabla\hat q(x_k)^*\lambda,
$$
and then consider the unique solution $w_k$ to the quadratic program \eqref{quad} with $q$ replaced by $\Hat q$. Using the same arguments as in the proof of Theorem~\ref{PropNecTiltStab2Reg} gives us the convergence $w_k\to w$ as $k\to\infty$ and the following relationships held for all $k\in\N$: $-w_k\in\hat K(x_k,x_k^\ast)$,
\[
\widehat N_{{\rm gph}\,\partial\delta_{\hat\Gamma}}(x_k,x_k^\ast)=\widehat N_{{\rm gph}\,\hat\partial\delta_{\hat\Gamma}}(x_k,x_k^\ast)=\big\{(u^\ast,u)\big|\;u\in\hat K(x_k,x_k^\ast),\;u^\ast\in-\nabla^2\big\la\lambda,q\big\ra(x_k)u+ \hat K(x_k,x_k^\ast)^*\big\},
\]
and $\nabla_x^2\hat\Lag(x_k,\lambda)w_k\in\breve\partial^2\hat f(x_k,\nabla\varphi(x_k)+x_k^\ast)(w_k)$, where
$$
\hat\Gamma:=\big\{x\in\R^n\big|\;\hat q(x)\in\Theta\},\;\hat f:=\varphi+\delta_{\hat\Gamma},\;\hat\Lag(\cdot,\lambda):=\varphi+\la\lambda,\hat q\ra,
$$
and $\hat K$ denotes the critical cone \eqref{crit-cone} generated by the aforementioned hat-constructions. Since we obviously have
$\nabla\varphi(x_k)+x_k^\ast=\nabla_x\hat\Lag(x_k,\lambda)\to 0$ as $k\to\infty$ as well as
\[
\lim_{k\to\infty}\big\la w_k,\nabla_x^2\hat\Lag(x_k,\lambda)w_k\big\ra=\big\la w,\nabla_x^2\hat\Lag(\xb,\lambda)w\big\ra\le 0,
\]
it follows from Theorem~\ref{ThTiltStab} that $\xb$ is not a tilt-stable minimizer of \eqref{EqPertProbl}, and we are done. $\h$

\section{Discussions and Examples}
\sce

In this section we discuss some remarkable features of the obtained second-order sufficient conditions and characterizations of tilt-stable minimizers in NLPs as well as the imposed MSCQ and BEPP qualification conditions, which ensure their validity. The presented examples reveal striking differences between the new results and those known in the literature and also illustrate new phenomena on tilt stability that have not been observed earlier.

Recall that the first {\em characterization} of tilt-stable minimizers in NLPs is obtained in \cite[Theorem~5.2]{MoRo12} {\em under LICQ} in the pointbased form of the classical {\em SSOSC} \cite{Rob80}:
\begin{eqnarray}\label{tilt-mr}
\big\la w,\nabla_x^2\Lag(\xb,\lambda)w\big\ra>0\;\mbox{ whenever }\;w\ne 0\;\mbox{ with }\;\la\nabla q_i(\xb),w\ra=0\;\mbox{ for all }\;i\in E\cup I^+(\lambda),
\end{eqnarray}
where $\lm\in\R^l$ is the {\em unique} Lagrange multiplier satisfying the KKT system \eqref{EqFO}. It has been well recognized that the simultaneous fulfillment of LICQ and SSOSC is a {\em characterization} of Robinson's {\em strong regularity} \cite{Rob80} for the variational inequality associated with KKT \eqref{EqFO}, and thus tilt stability of the local minimizer $\ox$ in \eqref{EqNonlProgr} is {\em equivalent} to strong regularity of $\ox$ in \eqref{EqFO} under the validity of LICQ, which is a necessary condition for strong regularity; see \cite[Corollary~5.3]{MoRo12} with the references and discussions therein. All the examples presented below demonstrate that in the results obtained in this paper in the absence of LICQ, which is surely not mandatory for tilt-stable minimizers, the property of tilt stability is {\em far removed} from strong regularity while postulating nevertheless a nice behavior of local minimizers from both qualitative and quantitative/numerical viewpoints.

It is shown in \cite[Theorem~3.5]{MoOut13} that SSOSC \eqref{tilt-mr}, assumed to hold for {\em all} the Lagrange multipliers in \eqref{EqFO}, is still a {\em sufficient} condition for tilt-stable minimizers in NLPs with inequality constraints when LICQ is relaxed to the simultaneous fulfillment of MFCQ and CRCQ at the reference point. The subsequent result of \cite[Theorem~4.3]{MoNg14} provides a characterization of tilt-stable minimizers in the same setting as in \cite{MoOut13} while being expressed via the {\em non-pointbased} USOSC discussed above in Section~1. Furthermore, \cite[Example~4.5]{MoNg14} demonstrates that the pointbased SSOSC {\em fails} in this setting, i.e., it cannot recognize a tilt-stable minimizer under MFCQ and CRCQ. It is worth mentioning to this end that the major difference of SSOSC \eqref{tilt-mr} from the similarly looking condition \eqref{2nd-nec3e} is that the positive-definiteness of the Hessian $\nabla_x^2\Lag(\xb,\lambda)$ in the latter one is required for the {\em larger} index set $I^+$ from \eqref{I++} {\em independent} of $\lm$. Thus condition \eqref{2nd-nec3e} is weaker than \eqref{tilt-mr} providing a {\em pointbased characterization} of tilt stability under the validity of CRCQ by Theorem~\ref{ThTiltStabCRCQ} while SSOSC fails to do it even under assuming in addition that MFCQ holds at this point.\vspace*{0.05in}

We now show that the usage of the new sufficient condition \eqref{2nd-suf3} from Corollary~\ref{ThTiltStabSuff1}, which involves not all the Lagrange multiplies but only the {\em extreme} ones in {\em critical directions} $\lm\in\bar\Lambda_{\cal E}$ from \eqref{2nd-suf2}, allows us to recognize a tilt-stable minimizer that does exist in \cite[Example~4.5]{MoNg14}.

\begin{Example}[\bf pointbased recognizing tilt stability via extreme multipliers in critical directions under MFCQ and CRCQ]\label{ExSSOSCwNotNec} {\rm Consider the following nonlinear program in $\R^3$:
\begin{eqnarray*}
\begin{array}{ll}
\mbox{minimize }&\varphi(x):=\disp\frac 14 x_1+x_3 +x_3^2-x_1x_2\;\mbox{ for }\;x=(x_1,x_2,x_3)\\
\mbox{subject to }&q_1(x):=x_1-x_3\le 0,\quad q_2(x):=-x_1-x_3\le 0,\\
&q_3(x):=x_2-x_3\le 0,\quad q_4(x):=-x_2-x_3\le 0.
\end{array}
\end{eqnarray*}
It is easy to check that MFCQ and CRCQ hold at $\xb=(0,0,0)$, and thus both MSCQ and BEPP are satisfied at $\bar x$ by Proposition~\ref{LemSuffCondBEPP}(ii). We can directly calculate that
\begin{eqnarray*}
&&\Lambda\big(\xb,-\nabla\varphi(\xb)\big)=\big\{\lambda\in\R^4_+\big|\;\lambda_1-\lambda_2=-\frac 14,\;\lambda_4=\lambda_3,\;\lambda_1+\lambda_2+\lambda_3+\lambda_4=1\big\},\\
&&\E\big(\xb,-\nabla\varphi(\xb)\big)=\Big\{\Big(0,\frac 14,\frac 38,\frac 38\Big),\Big(\frac 38,\frac 58,0,0\Big)\Big\},\;\mbox{ and }\;K\big(\ox,-\nabla\ph(\ox)\big)=\big(0,0,0\big).
\end{eqnarray*}
Hence the second-order sufficient condition \eqref{2nd-suf3} is trivially fulfilled due to $\bar\Lambda_{\cal E}=\emp$, and thus it recognizes tilt stability of the local minimizer $\ox$ in this example.}
\end{Example}

\begin{Remark} [\bf other consequences of Example~\ref{ExSSOSCwNotNec}]\label{ex8.1} {\rm Besides the main purpose of Example~\ref{ExSSOSCwNotNec}, it allows us to illustrate some other remarkable phenomena on tilt stability.

{\bf (i)} The tilt-stable minimizer $\ox$ in Example~\ref{ExSSOSCwNotNec} cannot be recognized not only by SSOSC \eqref{tilt-mr}, but also by its  {\em relaxed version} involving extreme multipliers:
\begin{eqnarray}\label{EqSSOSCw}
\begin{array}{ll}
\big\la w,\nabla^2_x\Lag(\xb,\lambda)w\big\ra>0\;\mbox{ whenever }\;\lambda\in\Lambda\big(\xb,-\nabla\varphi(\xb)\big)\cap\E\big(\xb,-\nabla\varphi(\xb)\big),\\
0\not=w\in\R^n,\;\mbox{ and }\;\big\la\nabla q_i(\xb),w\big\ra=0\;\mbox{ for all }\;i\in E\cup I^+(\lambda),
\end{array}
\end{eqnarray}
which differs from our new condition \eqref{2nd-suf3} by omitting the {\em critical directions} in the construction of $\bar\Lambda_{\cal E}$. Indeed,
taking $\lambda=(\frac 38,\frac 58,0,0)\in\Lambda(\xb,-\nabla\varphi(\xb))\cap\E(\xb,-\nabla\varphi(\xb))$ and $w=(0,1,0)$ in the setting of Example~\ref{ExSSOSCwNotNec}, we arrive at the relationships
\begin{eqnarray*}
\big\la\nabla q_i(\xb),w\big\ra=0\;\mbox{ for }\;i\in\{1,2\}=I^+(\lambda)\;\mbox{ while }\;\big\la w,\nabla_x^2\Lag(\xb,\lambda)w\big\ra=0,
\end{eqnarray*}
which show that the ``non-critical" counterpart \eqref{EqSSOSCw} of \eqref{2nd-suf3} fails at the tilt-stable minimizer $\xb$.

{\bf (ii)} Example~\ref{ExSSOSCwNotNec} cannot be directly used to illustrate Theorem~\ref{ThTiltStabCRCQ}, since the critical cone is trivial in this example while the opposite is assumed in the theorem. However, increasing the dimension of the problem by adding the term $\frac 12 x_4^2$ to the cost function in Example~\ref{ExSSOSCwNotNec} gives us an NLP with $K(\ox,-\nabla\ph(\ox))=(0,0,0)\times\R\ne\{0\}$ at the tilt-stable minimizer $\ox=0$ and such that the new condition \eqref{2nd-nec3e} holds while SSOSC \eqref{tilt-mr} fails therein. Indeed, in this case we have $I^+=\{1,2,3,4\}$, and therefore $[\la\nabla q_i(\ox),w\ra=0\;\mbox{ for all }\;i\in I^+]$ implies that $w=(0,0,0,w_4)$ and $\big\la w,\nabla^2_x\Lag(\xb,\lambda)w\big\ra=w_4^2$ for all $\lambda\in\Lambda\big(\ox,-\nabla\ph(\ox)\big)$, i.e.,  condition \eqref{2nd-nec3e} is satisfied. On the other hand, we get $I^+(3/8,5/8,0,0)=\{1,2\}$, which shows the violation of \eqref{tilt-mr} for $w=(0,1,0,0)$.}
\end{Remark}

The next example reveals the situation when both MFCQ and CRCQ fail at a local minimizer $\ox$ while SOSCMS \eqref{SOSCMS}, and hence MSCQ and BEPP by Theorem~\ref{ThSSOSCMS}, are satisfied at this point together with the other assumptions of Theorem~\ref{CorEquivTiltStab2reg} ensuring therefore that the second-order condition \eqref{2nd-suf3} provides a complete {\em pointbased characterization} of tilt stability for $\ox$.

\begin{Example}[\bf pointbased characterization of tilt stability under 2-regularity but without MFCQ and CRCQ]\label{ExSuffNec} {\rm Given a parameter pair $(a,b)\in\R^2$, consider the following NLP in $\R^3$:
\begin{eqnarray}\label{EqProblNonDeg}
\begin{array}{ll}
\mbox{minimize }&\varphi(x):=-x_1+\disp\frac a2 x_2^2+\frac b2x_3^2\;\mbox{ for }\;x=(x_1,x_2,x_3)\\
\mbox{subject to }&q_1(x):=x_1-\disp\frac 12 x_2^2\le 0,\quad q_2(x):=x_1-\frac 12 x_3^2\le 0,\\
&q_3(x):=-x_1-\disp\frac 12 x_2^2-\frac 12 x_3^2\le 0.
\end{array}
\end{eqnarray}
Letting $\xb=(0,0,0)$, it is easy to observe that both MFCQ and CRCQ are violated at $\ox$ while SOSCMS \eqref{SOSCMS} holds with $E_2=I_2=\emp$ therein. To check the latter, pick any vectors $0\not=u=(u_1,u_2,u_3)$ and $0\not=(\lambda_1,\lambda_2,\lambda_3)\in N_\Theta(q(\xb))$ satisfying $\nabla q(\xb)u\in T_\Theta(q(\xb))$ and $\la\nabla q(\xb),\lambda\ra=0$ and then get $u_1=0$, $\lambda_1+\lambda_2-\lambda_3=0$, $\lambda_1\ge 0$, $\lambda_2\ge 0$, $\lambda_3\ge 0$, and so $\lambda_3=\lambda_1+\lambda_2>0$. This gives us
\[
\big\la u,\nabla^2\la\lambda,q\ra(\xb)u\big\ra=-\big(\lambda_1+\lambda_3\big)u_2^2-\big(\lambda_2+\lambda_3\big)u_3^2<0
\]
and thus verifies the validity of SOSCMS in this setting.

The corresponding set of multipliers \eqref{Lambda} and its extreme points are calculated by, respectively,
$$
\Lambda\big(\xb,-\nabla\varphi(\xb)\big)=\big\{\lambda\in\R^3_+\big|\;\lambda_1+\lambda_2-\lambda_3=1\big\},\quad\E\big(\xb,-\nabla\varphi(\xb)\big)=
\big\{(1,0,0),(0,1,0)\big\}.
$$
The critical cone amounts to $K(\xb,-\nabla\varphi(\xb))=\{0\}\times\R\times\R$, and for $0\not=v\in K(\xb,-\nabla\varphi(\xb))$ we have
\[
\Lambda\big(\xb,-\nabla\varphi(\xb);v\big)=\begin{cases}(1,0,0)&\mbox{if $v_2^2<v_3^2$,}\\
\big\{(\lambda_1,\lambda_2,0)\in\R^3_+\big|\;\lambda_1+\lambda_2=1\big\}&\mbox{if $v_2^2=v_3^2$,}\\
\big(0,1,0\big)&\mbox{if $v_2^2>v_3^2$.}
\end{cases}
\]
This tells us that $\Lambda(\xb,-\nabla\varphi(\xb);v)$ is a singleton when $v_2^2\not=v_3^2$, and thus we meet the assumptions of Theorem~\ref{PropNecTiltStab2Reg}(b) used also in Theorem~\ref{CorEquivTiltStab2reg} by showing that for every $0\ne v\in K(\xb,-\nabla\varphi(\xb))$ with $v_2^2=v_3^2$ and every $\lambda\in\Lambda_{\cal E}(\xb,-\nabla\varphi(\xb);v)=\{(1,0,0),(0,1,0)\}$ there is a maximal subset $\hat\J$ of $\J(\xb;v)$ such that $I^+(\lambda)\subset\hat\J$ and $(q_i)_{i\in E\cup\hat\J}$ is 2-regular at $\ox$ in the direction $v$. To proceed, observe from the above that $v_1=0$ and $v_2^2=v_3^2\not=0$ in our case and that the set $\Xi$ from \eqref{Xi} is
\[
\Xi(\xb;v)=\big\{(z_1,z_2,z_3)\big|\;z_1-v_2^2\le 0,\;z_1-v_3^2\le 0,\;-z_1-v_2^2-v_3^2\le 0\big\},
\]
which gives us $\J(\xb;v)=\{\emp,\{3\},\{1,2\}\}$. Then we have that $\hat\J=\{1,2\}$ is a maximal element of $\J(\xb;v)$,
$I^+(\lambda)\subset\hat\J$, $\lambda\in\{(1,0,0),(0,1,0)\}$, and for every $\alpha\in\R^2$ the system
\begin{eqnarray*}
&&\big\la\nabla q_1(\xb),u\big\ra+\big\la v,\nabla^2q_1(\xb)w\big\ra=u_1-v_2w_2=\alpha_1,\;\big\la\nabla q_1(\xb),w\big\ra=w_1=0,\\
&&\big\la\nabla q_2(\xb),u\big\ra+\big\la v,\nabla^2q_2(\xb)w\big\ra=u_1-v_3w_3=\alpha_2,\;\big\la\nabla q_2(\xb),w\big\ra=w_1=0
\end{eqnarray*}
has a solution $(u,w)$, e.g., $u=(\alpha_1,0,0)$ and $w=(0,0,(\alpha_1-\alpha_2)/v_3)$. This verifies the required 2-regularity in Theorem~\ref{CorEquivTiltStab2reg}, and so we can apply the tilt-stability characterizations therein. The straightforward second-order calculation in the positive-definiteness condition \eqref{2nd-suf3} shows that $\xb$ is a tilt-stable local minimizer in \eqref{EqProblNonDeg} {\em if and only if} $a>1$ and $b>1$. Furthermore, we can compute the exact bound of tilt stability of $\ox$ in this program by tilt$(\ph,q,\ox)=1/\min\{a-1,b-1\}$.

Note finally that in this example SSOSC \eqref{tilt-mr} fails at $\ox$ if $a=b=2$, $\lambda=(0,2,1)\in\Lambda(\xb,-\nabla\varphi(\xb))$, and $w=(0,0,1)$. Indeed, we have then $\la\nabla q_2(\xb),w\ra=\la\nabla q_3(\xb),w\ra=0$ while $\la w,\nabla^2_x\Lag(\xb,\lambda)w\ra=-1$.}
\end{Example}

The next example demonstrates that the additional assumptions of Theorem~\ref{CorEquivTiltStab2reg} (taken from Theorem~\ref{PropNecTiltStab2Reg}) ensuring the {\em necessity} of the second-order sufficient condition  \eqref{2nd-suf3} for tilt-stable minimizers, {\em cannot be dropped} even under the validity of MFCQ.

\begin{Example}[\bf nondegeneracy and 2-regularity are essential for pointbased characterizing tilt-stable minimizers]\label{ExNotNec1}{\rm
Consider the the following NLP in $\R^3$:
\begin{eqnarray}\label{exa3}
\begin{array}{ll}
\mbox{minimize }&\ph(x):=-x_1+\disp\frac 12 x_2^2,\quad x=(x_1,x_2,x_3),\\
\mbox{subject to }&q_1(x):=x_1+x_3^2\le 0,\;q_2(x):=x_1\le 0.
\end{array}
\end{eqnarray}
We obviously have that {\em MFCQ holds} at $\ox=(0,0,0)$, and hence both MSCQ and BEPP assumed in Theorem~\ref{ThTiltStabSuff} are satisfied at this point. Since the second constraint in \eqref{exa3} is clearly redundant, we can consider the equivalent version of this problem {\em without} the latter constraint and easily deduce from Theorem~\ref{ThTiltStabSuff} that $\ox$ is a tilt-stable minimizer in it with modulus $\kk=1$. However, applying Theorem~\ref{ThTiltStabSuff} to the original (``full") version of \eqref{exa3} shows that the second-order sufficient condition \eqref{2nd-suf3} fails, and so we cannot make a conclusion about tilt stability of $\ox$ in \eqref{exa3} by using this theorem. Indeed, taking $v=(0,1,0)$,  $\lambda=(0,1)$, and $w=(0,0,1)$ gives us
\begin{eqnarray*}
\begin{array}{ll}
v\in K\big(\xb,-\nabla\varphi(\xb)\big)=\big\{v\in\R^3\big|\;v_1=0\big\},\;\lambda\in\Lambda_{\cal E}\big(\xb,-\nabla\varphi(\xb);v\big)=\big\{(1,0),(0,1)\big\},\\
I^+(\lambda)=\big\{2\big\},\;\big\la\nabla q_2(\xb),w\big\ra=0,\;\mbox{ and }\;\big\la w,\nabla_x^2\Lag(\xb,\lambda)w\big\ra=0,
\end{array}
\end{eqnarray*}
which shows that the sufficient condition \eqref{2nd-suf3} for tilt stability is not fulfilled at $\ox$. The reason is that the additional assumption of Theorem~\ref{PropNecTiltStab2Reg} ensuring the necessity of \eqref{2nd-suf3} for tilt stability are not satisfied here. To see this, observe that
the set of Lagrange multipliers \eqref{crit-mult} in the critical direction $v$ is {\em not a singleton} (i.e., $\ox$ degenerates in this direction), which violates the assumption in Theorem~\ref{PropNecTiltStab2Reg}(a). Furthermore, the set of active inequality constraint indexes \eqref{active1} in this direction is $\J(\xb;v)=\{\emp,\{1,2\}\}$, which shows that the {\em 2-regularity} assumption of Theorem~\ref{PropNecTiltStab2Reg}(b) is also violated.}
\end{Example}

The last example in this section is a modification of Example~\ref{ExNotNec1} illustrating the phenomenon on tilt stability revealed in Theorem~\ref{ThNecCond2}, which shows that there are two NLPs with the same cost function and the same values of the constraints functions and their derivatives up to the second order at the reference point such that this point satisfies SOSCMS and gives a tilt-stable local minimum for one program but not for the other one. The example presented below illustrates this phenomenon under MFCQ (which is stronger than SOSCMS) in the case where the constraint functions and their derivatives up the third order are the same in at the point in question. Actually this example can be further modified to exhibit the aforementioned phenomenon under the validity of MFCQ in the case where the constraint functions and their derivative of {\em any order} are the same at the reference point.

\begin{Example}[\bf pointbased characterizations of tilt stability are not possible under MFCQ alone]\label{ExNotNec2} {\rm Consider the following NLP in $\R^3$, which differs from \eqref{exa3} by the term $-x^4_2$ in the function $q_1(x)$:
\begin{eqnarray}\label{exa4}
\begin{array}{ll}
\mbox{minimize }&\ph(x):=-x_1+\disp\frac 12 x_2^2,\quad x=(x_1,x_2,x_3),\\
\mbox{subject to }&q_1(x):=x_1-x^4_2+x_3^2\le 0,\;q_2(x):=x_1\le 0.
\end{array}
\end{eqnarray}
We obviously have that MFCQ holds at $\ox=(0,0,0)$ in \eqref{exa4} and the values of the constraint functions and their derivatives up to the third order at $\ox$ are same in \eqref{exa3} and \eqref{exa4}. As shown in Example~\ref{ExNotNec1}, $\ox$ is a tilt stable minimizer of \eqref{exa3} while the second-order sufficient condition \eqref{2nd-suf3} fails for $\ox$ in \eqref{exa3}.

To verify that $\ox$ is {\em not} a tilt-stable minimizer for NLP in \eqref{exa4}, pick the same elements  $v=(0,1,0)$, $\lambda=(0,1)$, and $w=(0,0,1)$ as in Example~\ref{ExNotNec1} and then, according to Definition~\ref{DefTiltStab} of tilt stability and its adjustment for NLPs in Section~6, consider the problem
\begin{eqnarray*}
\mbox{minimize }\;-x_1+\disp\frac 12 x_2^2-ux_2\;\mbox{ subject to }\;x_1-x^4_2+x_3^2\le 0,\;x_1\le 0
\end{eqnarray*}
with only one tilt parameter $u\in\R$ in this case. For each $u\ne 0$ the latter parametric optimization problem has two distinct solutions
$(0,u,\pm u^2)$, which excludes the validity of tilt stability of $\ox$ in \eqref{exa4}.

Note that we can also construct an NLP equivalent (in the second-order) to \eqref{exa3} but without tilt stability at $\ox$ by using the proof of Theorem~\ref{ThNecCond2}. Indeed, let $z=v=(0,1,0)$ and thus get the functions
$$
\hat q_1(x)=x_1+x_3^2-\big(x_1^2+x_2^2+x_3^2\big)^2\;\mbox{ and }\;\hat q_2(x)=x_1
$$
in the proof therein, where $\hat q_1$ is surely more complicated in comparison with \eqref{exa4}.}
\end{Example}

\begin{Remark}[\bf tilt stability and critical multipliers]\label{cr-mult} {\rm Finally in this section, we discuss some {\em numerical consequences} of the obtained results on tilt stability. This concerns relationships between tilt stability of local minimizers in NLPs and the so-called {\em critical minimizers} that have been recently discovered and then strongly investigated in the excellent book by Izmailov and Solodov \cite{IS14}. It is shown in \cite{IS14} that critical multipliers, which may appear even in the case of unique multipliers under LICQ, are largely responsible for {\em slow convergence} of major primal-dual numerical algorithms including Newton and Newton-type methods, the augmented Lagrangian method, the sequential quadratic programming method, etc. Therefore it is highly desired from the numerical viewpoint to {\em rule out} the existence of critical multipliers and so to be able making such a conclusion based on the {\em initial data} of the NLP in question. These and related issues have been discussed in the recent comments of the second author \cite{m15} on the survey by Izmailov and Solodov devoted to critical multipliers, which is based on their book \cite{IS14}. It is {\em conjectured} in \cite{m15} that under appropriate qualification conditions tilt stability excludes the existence of critical multipliers.

The results obtained in this paper shed light on this conjecture and its consequences for primal-dual algorithms of numerical optimization. Indeed, it can be derived from \cite[Theorem~4.3]{MoNg14} that the simultaneous validity of MFCQ and CRCQ (and surely in the case of LICQ) at the given tilt minimizer $\ox$ ensures that critical multipliers {\em do not appear} at $\ox$, i.e., the above conjecture is valid in this setting. Thus the pointbased {\em necessary} conditions for (as well as the characterizations of) tilt-stable minimizers established in Section~7 allow us to exclude, under the validity of MFCQ and CRCQ at $\ox$, undesired behavior of the aforementioned numerical algorithms. Observe, in particular, that our major pointbased second-order condition \eqref{2nd-suf3}, which characterizes the tilt stability of $\ox$ by Theorem~\ref{CorEquivTiltStab2reg}(ii) via the positive-definiteness of the Hessian of the Lagrange function only for {\em extreme}  multipliers in {\em critical directions} \eqref{2nd-suf2}, tells us now that {\em all} the Lagrange multipliers are {\em noncritical} at $\ox$ in the sense of \cite{IS14} in this rather general setting.

On the other hand, {\em MFCQ alone} does not allow us to exclude the existence of critical multipliers at tilt stable minimizers. It happens, in particular, in the setting of Example~\ref{ExNotNec1} under MFCQ and also in Example~\ref{ExSuffNec} under the weaker SOSCMS. Alexey Izmailov (private communication) informed us about a two-dimensional example admitting the {\em unique} critical Lagrange multiplier under the validity of MFCQ (but not LICQ) at a tilt-stable minimizer. Thus the question remains on what (weaker than CRCQ) should be added to MFCQ, or even what can replace MFCQ and CRCQ together, to ensure that tilt stability excludes critical multipliers at the reference local minimizer.}
\end{Remark}

\section{Open Questions and Further Research}

It seems to us that this paper basically clarifies the situation with second-order necessary and sufficient conditions for tilt-stable local minimizers in finite-dimensional NLPs, and not much is expected to be added to this theory. However, principal questions remain about using the obtained results and the very notion of tilt stability in {\em numerical optimization} including, in particular, more work on relationships between tilt stability and critical multipliers discussed at the end of Section~8. Challenging issues arise on infinite-dimensional (mainly Hilbert space) extensions of the obtained pointbased characterizations and also on establishing appropriate counterparts of the NLP tilt stability theory above in other remarkable classes of constrained optimization, particularly for problems of {\em conic programming}.

Among the most natural topics of the future research we mention developing a comprehensive theory of {\em full stability} for local minimizers in NLPs as well as in other classes of constrained optimization and variational problems. The notion of full stability was introduced by Levy, Poliquin and Rockafellar \cite{LPR00} in the extended-real-valued framework of unconstrained optimization as a far-going generalization of tilt stability. Recently it has been largely extended to various classes of constrained optimization problems in \cite{MN14,mnr14,mos14,mrs13,ms14}. However, most of the results obtained in these papers impose {\em nondegeneracy} assumptions (analogs of LICQ) on the corresponding constraints. The only exception is \cite{MN14}, where neighborhood characterizations of full stability in NLPs are obtained under the simultaneous validity of partial versions of MFCQ and CRCQ. A major goal of the future research is to extend the theory of tilt stability developed in this paper to the case of fully stable local minimizers in NLPs. Note that full stable minimizers seem to be more appropriate to rule out critical multipliers according to the second conjecture in \cite{m15}.

\section*{Acknowledgements}
The research of the first author was supported by the Austrian Science Fund (FWF) under grant P26132-N25. The research of the second author was partially supported by the USA National Science Foundation under grants DMS-12092508 and DMS-1512846 and by the USA Air Force Office of Scientific Research under grant No.\,15RT0462.

\end{document}